\documentclass{amsart}
\usepackage{graphicx}
\usepackage{natbib}

\newtheorem{thm}{Theorem}[section]
\newtheorem{cor}[thm]{Corollary}
\newtheorem{lem}[thm]{Lemma}

\newtheorem{defn}[thm]{Definition}

\newcommand{\diff}{\mathrm{d}}

\def\M{\mathcal M}
\begin{document}

\title[]{Empirical geodesic graphs and CAT(k) metrics for data analysis}%
\author{Kei Kobayashi}%
\address{Keio University}%
\email{kei@math.keio.ac.jp}%

\author{Henry P Wynn}%
\address{London School of Economics}%
\email{h.wynn@lse.ac.uk}%

\keywords{intrinsic mean, extrinsic mean, CAT(0), curvature, metric cone, cluster analysis, non-parametric analysis}%

\begin{abstract}
A methodology is developed for data analysis based on empirically constructed geodesic metric spaces.
For a probability distribution, the length along a path between two points can be defined as the
amount of probability mass accumulated along the path. The geodesic, then, is the shortest such path and
defines a geodesic metric. Such metrics are transformed in a number of ways to produce parametrised families of geodesic metric spaces,  empirical versions
of which allow computation of intrinsic means and associated measures of dispersion. These reveal properties of the data, based on geometry, such as those that are difficult to see from the raw Euclidean distances. Examples of application include clustering and classification.
For certain parameter ranges, the spaces become CAT(0) spaces and the intrinsic means are unique. In one case, a minimal spanning tree of a graph based on the data becomes CAT(0). In another, a so-called ``metric cone" construction allows extension to CAT($k$) spaces. It is shown how to empirically tune the parameters of the metrics, making it possible to apply them to a number of real cases.

\end{abstract}
\maketitle
\begin{center}
{\bf This paper is to appear in \textit{\textbf{Statistics and Computing}}, 2019, \\
DOI 10.1007/s11222-019-09855-3.}
\end{center}

\section{Introduction}

In much statistics and data analysis, the metric (distance) for data points is fixed and the loss function is selected from a set of candidates for loss functions and/or tuned by a parameter.   
However, in the paper, we fix a loss function (usually the squared loss) and instead select/tune the metric. 
The motivation of such metric-based approach is to propose a
set of metrics and a method to select a metric from it which are naturally acquired from geometrical aspects. 
This enables us to import huge existing literature of various branches of geometry into data analysis.
We will begin by focusing the curvature.
In this section, after explanation by a motivational example, exiting studies related to our geometrical approach will be surveyed.

\subsection{Example}

\begin{figure}[tbp]
\begin{center}
 \begin{minipage}{4.1cm}
  \begin{center}
   \includegraphics[width=4.4cm]{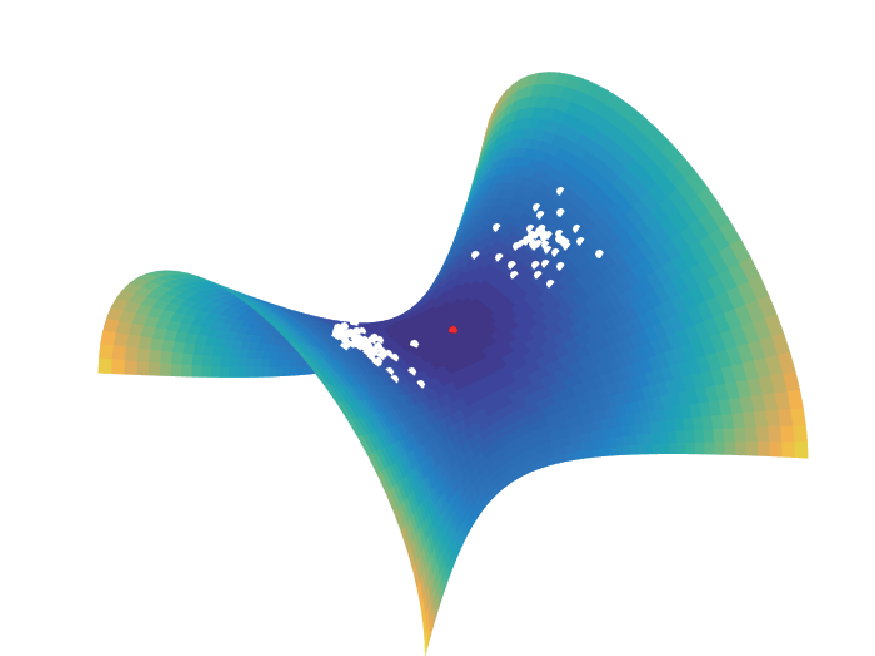}\\
   (a) Hyperboloid
  \end{center}
 \end{minipage}
 \begin{minipage}{4.1cm}
  \begin{center}
   \includegraphics[width=4.4cm]{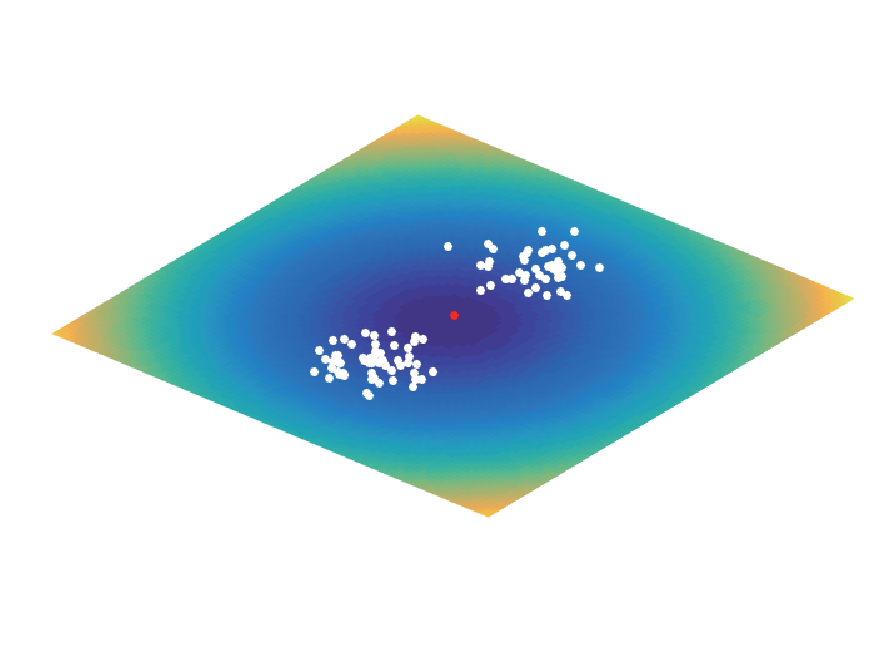}\\
   (b) Plane
  \end{center}
 \end{minipage}
 \begin{minipage}{4.1cm}
  \begin{center}
   \includegraphics[width=4.4cm]{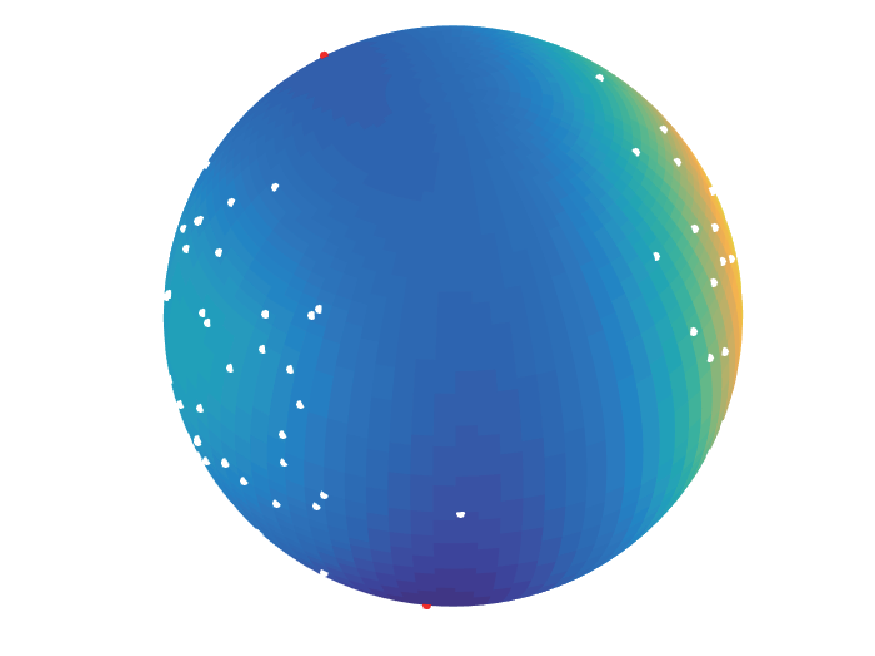}\\
   (c) Sphere
  \end{center}
 \end{minipage}
\caption{The Fr\'echet function $f(m)$ for data (white dots) on (a) a hyperboloid (curvature $c=-1$),
(b) a plane ($c=0$) and (c) a sphere ($c=1$). The bluer represents the smaller value of $f(m)$. 
The red dots represent the local minima of $f(m)$. Only for the sphere, $f(m)$ has multiple minima.
}
\label{fig:hyperboloid}
\end{center}
\end{figure}

For a random variable $X$ on a metric space $\M$ endowed with a  metric $d( \cdot, \cdot)$ the general  {\em intrinsic
mean} is defined by
\begin{equation*}
\mu = \mbox{arg} \min_{m \in \M} \mbox{E}[d(X,m)^2].
\end{equation*}
The empirical intrinsic mean based on data $x=\{x_1, \ldots, x_n\}$, sometimes called the Fr\'echet mean, is defined as
\begin{equation*} \hat{\mu} = \mbox{arg} \min_{m \in \M} f (m),\end{equation*}
where
\begin{equation*}f(m) = \sum_{i=1}^n d(x_i,m)^2.\end{equation*}
The function $f(m)$ is sometimes referred to as the Fr\'echet function.  
For Euclidean space, $\hat{\mu} = \bar{x}$, the sample mean. 
In general, $f(m)$ is not necessarily convex and the means, $\hat{\mu}$, are not unique.
Figure \ref{fig:hyperboloid} shows that the curvature can affect the property of $f(m)$.
In particular, for so-called CAT(0) spaces, which (trivially) include Euclidean spaces, the intrinsic means $\hat{\mu}$ are unique.

Even when the mean is not unique, the function $f(m)$ can yield useful information, for example about clustering. We can also define second-order quantities:
\begin{equation*}s_0^2 = \inf_{m \in \M}  \frac{1}{n}\sum_{i=1}^n d(x_i,m)^2 = \frac{1}{n}\sum_{i=1}^n d(x_i,\hat{\mu})^2\end{equation*}
and
\begin{equation*}s_{1}^2 = \frac{2}{n(n-1)} \sum_{i<j}^n d(x_i,x_j)^2.\end{equation*}
The quantity $s_0^2$ is sometimes called the  Fr\'echet variance. We name $s_1^2$ as the mean pairwise discrepancy.

A key concept in the study of these issues is that the metrics are global geodesic metrics, that is metrics based on
the shortest path between points measured by integration along a path with respect to a local metric.
The interplay between the global and the local will concern us to a considerable extent.


The general form of the Fr\'echet function depends, here, on three parameters, $\alpha, \beta, \gamma,$ and it can be written in compact form:
\begin{equation*}f_{\alpha, \beta, \gamma}(m) = \sum_{i=1}^n \{ g_\beta(d_{\alpha}(x_i,m))\}^{\gamma},\end{equation*}
where the function $g_{\beta}$ and the construction of $d_{\alpha}$ are given below.
Once we have introduced this new class of metrics, variety of statistics can be generalised:
intrinsic mean, variance, clustering (based on local minima of $f(m)$).
For classification problems, we can select an appropriate metric by cross-validation.

There are many ways to transform one metric into another, 
regardless of whether they are geodesic metrics.  A straightforward way is to use a concave function
$g$ such that given a metric $d(\cdot, \cdot)$, the new metric is $d'(\cdot,\cdot) = g(d(\cdot, \cdot))$. This is plausible if we use non-convex $f_{\alpha,\beta, \gamma}$, which are useful, as will be explained, in clustering and classification. 
Such concave maps are often interpreted as loss functions, but we will consider
them in terms of changes of metric which may lead to selection using geometric concepts. This is particularly true for the construction based on the $g_{\beta}$ in Section \ref{sec:alpha} of the paper.
In Table \ref{table:generalised_statistics}, we summarise such generalised statistics.
\begin{table*}[thb]
\begin{center}
\caption{A summary of generalised statistics by introducing $\alpha,\beta$ and $\gamma$.}
\label{table:generalised_statistics}
\begin{tabular}{|c||c|c|}
\hline
& Euclidean
& Generalised metric\rule[-2mm]{0mm}{6mm}\\
\hline
Metrics
& {\small $d(x,y)=\|x-y\|$}
& {\small $d_{\alpha\beta}(x,y) = g_\beta(d_\alpha(x,y))$}\rule[-4mm]{0mm}{10mm}\\
\hline
Intrinsic mean
& {\small $\displaystyle\arg\min_{m \in \mathbb{E}^d} \sum_{i=1}^n \|x_i-m\|^2$}
& {\small $\displaystyle\arg\min_{m \in \mathcal{M}} \sum_{i=1}^n g_\beta(d_\alpha(x_i,m))^\gamma$}\rule[-5mm]{0mm}{12mm}\\
\hline
Variance
& {\small $\displaystyle\min_{m \in \mathbb{E}^d} \frac{1}{n} \sum_{i=1}^n \|x_i-m\|^2$}
& {\small $\displaystyle\min_{m \in \mathcal{M}} \frac{1}{n} \sum_{i=1}^n g_\beta(d_\alpha(x_i,m))^\gamma$}\rule[-5mm]{0mm}{12mm}\\
\hline
Fr\'echet function
& {\small $f(m)=\displaystyle\sum_{i=1}^n \|x_i-m\|^2$}
& {\small $f_{\alpha\beta\gamma}(m)=\displaystyle\sum_{i=1}^n g_\beta(d_\alpha(x_i,m))^\gamma$}\rule[-5mm]{0mm}{12mm}\\
\hline
\end{tabular}
\end{center}
\end{table*}

The basic definition and construction from a geodesic metric space to the special geodesics based on accumulation of density are given in the next section, together with the definition of a CAT(0) space.
In Section \ref{sec:geodesics}, we first show that means and medians in simple one-dimensional statistics can be placed into our framework.
Because geodesics themselves are one-dimensional paths, this should provide some essential motivation.
The $d_\alpha$-metric is obtained by a local dilation.
Our computational shortcut is to use {\em empirical} graphs, whose vertices are data points.

We will need, therefore, to define empirical geodesics. We start with a natural
geodesic defined via a probability density function in which the distance along a path is the amount of density ``accumulated" along that
path. Then, an empirical version is defined whenever a density is estimated. 

In Section \ref{sec:beta}, the $d_{\beta}$ metric is introduced.
It is based on a function derived from a geodesic metric via shrinking, pointwise, to an abstract origin (apex); that is to say an abstract cone is attached.
The smaller the value of $\beta$, the closer to the origin.
We cover the more general CAT($k$) spaces, giving some new results related to ``diameter" , in Section \ref{sec:cat_k}, including conditions for the uniqueness of intrinsic means not requiring the spaces to be CAT(0).

Section \ref{sec:choosing} provides a summary of the effect of changing $\alpha$ and $\beta$. After some discussion of the selection of $\alpha$ and $\beta$ in Section \ref{sec:examples}, Section \ref{sec:conclusion} covers some examples.

\subsection{Related existing studies}\label{sec:survey}

Manifold learning is a group of nonlinear dimension reduction techniques
including well-studied methods such as Isomap (\cite{Tenenbaum2000-to}), Locally Linear Embedding (LLE) (\cite{Saul2003-ze}) and Laplacian Eigenmaps (\cite{Belkin2002-nf}).
Most manifold learning methods are based on the ``manifold hypothesis,'' which is an assumption that the data is distributed around a smooth manifold with a lower dimension embedded in a higher-dimensional vector space (usually Euclidean space).
There are some similarity between the methods proposed in this paper and manifold learning methods though the original motivation of the research is different;
both methods focus on the geometrical structure of an embedded data space and, furthermore, use the geodesic length (shortest path length) in an empirical graph as the distance between data points. 
Our methods have significant differences from the manifold learning methods.
First, we control the curvature of the data space for data analysis via changing the metric while the metric in manifold learning context is fixed and to be estimated. 
Second, sometimes more positively (or negatively) curved data space is preferable in contrast to a situation in most manifold learning methods which attempt to estimate the manifold by making the approximated empirical graph locally flat (Euclidean) as possible.
For more details of the manifold learning methods, there are good surveys, e.g. \cite{Yang2006-pi}, \cite{Cayton2005-xr} (with other metric learning methods such as kernel learning) and \cite{Bengio2013-xp}(with various other data representation for machine learning).

Statistical shape analysis (\cite{Kendall2009-lv}, \cite{Ramsay2007-mt}, \cite{Srivastava2016-uk}), also known as object oriented data analysis (\cite{Marron2014-kt}), has a long history after a pioneering work by \cite{Kendall1984-tr} on random segmentations.
Statistical shape analysis studies geometrical structure of the set (shape space) of possible populations which themselves have some particular shapes.
For various kinds of shape spaces, computation of center points as mean and median and statistical methods as PCA and Bootstrap tests have been studied (see, e.g. \cite{Dryden2016-st}).
In particular, Tree space for analyzing phylogenetic trees (\cite{Billera2001-ve}, \cite{Wang2007-cv}) is closely related to our research.
Tree space is a set of tree graphs and the set is embedded in a Euclidean space with a tacitly defined metric.
The space is proved to have the CAT(0) property and therefore both geodesics between any pair of points and Fr\'echet mean of any finite data sets exit uniquely. 
Furthermore, a polynomial algorithm to compute the geodesics 
(\cite{Owen2011-jy}) and PCA on Tree space (\cite{Nye2011-tl}) have been proposed.
Besides such non-smooth spaces as Tree space, there have been many studies of non-parametric statistics on smooth manifolds. There are excellent textbooks in this area such as \cite{Bhattacharya2012-hv}, \cite{Patrangenaru2015-hg}.

Data analysis using Wasserstein distance also focuses geodesic distance  in the space of probability measures (see, e.g., \cite{Vallender1974-hs}, \cite{Villani2008-dd}, \cite{Peyre2018-jp}). 
Wasserstein Fr\'echet mean of measures is also studied (\cite{Cuturi2014-rr}) and uniqueness and computation of the mean depends on the curvature of the Wasserstein space.
For example, the 2-Wasserstein space for Gaussian measures has positive curvature in general (\cite{Takatsu2011-gj}) and therefore computation of Wasserstein Fr\'echet mean is difficult.
\cite{Panaretos2018-ea} is a useful survey of Wasserstein metric from statistical aspects and geometry of Wasserstein space is summarized in a section.
Because most of the recent algorithms in machine learning and computer graphics are based on some gradient methods, Wasserstein metric is becoming an active topic in such areas, e.g. Wasserstein-GAN (\cite{Arjovsky2017-yx}) and optimal transport of graphics with a penalty on the entropy (\cite{Solomon2015-ri}). 

Another area in statistics directly dealing with the curvatures is information geometry (\cite{Amari1985-sk}, \cite{McCullagh1986-zl}).
Fisher metric, based on the Fisher information matrix, is induced in a statistical model manifold which is curved in an embedding space (usually flat, for example, the space of an exponential family).
The asymptotic property and efficiency of estimators and predictors can be represented by the embedding curvature and naturally induced dual affine connections.
The role of curvatures in information geometry can be somewhat negative; even if the model manifold has non-zero embedding curvature, still some asymptotic efficiency of estimators (like bias-corrected MLE) can be proved.


As we have explained, there are many studies on statistics and data analysis using geodesics and curvatures, but our methods have some special features:
\begin{itemize}
\item The curvature of the data space holds not only by its own nature, but is controlled for data analysis.
\item The structure of empirical graphs used for computing the distance between the data points is not fixed but transformed via controlling the curvature.
\item Our methods can produce non-geodesic distances for data analysis from the aspects of curvatures though the curvatures cannot be defined for non-geodesic metric spaces.
This was achieved by considering the curvature of a metric space embedding the data space.
\end{itemize}

\section{Geodesics, intrinsic mean and extrinsic mean}\label{sec:geodesics}
The fundamental object in this paper is a geodesic metric space.
This is defined in two stages. First, define a metric space $\mathcal M=(X,d) $ with base space $X$ and metric $d(x,x')$. Sometimes, $\mathcal M$ will be a Euclidean space $E_d$ of dimension $d$,
containing the data points, but it may also be some special object such as a graph or manifold.
Second, define the length of a (rectilinear) path between two points $x,x' \in X$ and the
geodesic connecting $x$ and $x'$ as the shortest such path. 
The minimal length defines a metric $d^*(x,x')$, and the space  endowed with the geodesic metric is called the geodesic metric space, $\mathcal M^* = {\mathcal M}(X,d^*)$.

The interplay between $\mathcal M =  (X,d)$ and $\mathcal M^* =(X,d^*)$ will be critical for this paper, and,
as mentioned, we will have a number of ways of constructing $d^*$.

For data points $x_1, \ldots, x_n$ in $X$, the empirical intrinsic (Fr\'echet) mean is
\begin{equation*}\mu=\arg\inf_{\mu\in X} \sum_{i=1}^n d^*(x_i,\mu)^2.\end{equation*}
There are occasions when ${\mathcal M}^*$ can be represented as  a sub-manifold of a larger space (such as Euclidean space) ${\mathcal M}^+ = (X^+, d^+)$ with its own  metric $d^+$.
We can then talk about the extrinsic mean:
\begin{equation*}\mu^+=\arg\inf_{\mu \in X} \sum_{i=1}^n d^+(x_i,\mu)^2.\end{equation*}
Typically, the extrinsic mean is used as an alternative when the geodesic distance $d^*$ is hard to compute. The difficultly in considering the intrinsic mean in $X^+$ is that it may not lie in the original base space $X$. This leads to a third possibility, which is to
project it  back to $X$, in some way, as an approximation to the intrinsic mean $\mu$ (which may be hard to
compute). We will discuss this again in Section 4. See \cite{bhattacharya-2012} for further discussion on the intrinsic and extrinsic means.

\subsection{CAT(0) and CAT($k$) spaces}\label{cat}
CAT(0) spaces, which correspond to non-positive curvature Riemannian spaces, are important here because their intrinsic means are unique. The CAT(0) property is as follows.
Take any three points $\{a,b,c\}$ in a geodesic metric space $X$ and consider the ``geodesic triangle" of the points based
on the geodesic segments connecting them.
Construct a triangle in Euclidean 2-space with vertices $\{a', b' ,c'\}$, called the comparison triangle, whose Euclidean distances,
$\|a'-b'\|, \|b'-c'\|, \|a'-c'\|$, are the same as the corresponding geodesic distances just described: $d(a,b) = \|a'-b'\|$, etc. On the geodesic triangle select a point $x$ on the geodesic edge
between $b$ and $c$ and find the point $x'$ on the edge $b'c'$
of the Euclidean triangle such that $d(b,x) = \|b'-x'\|$.
Then the CAT(0) condition is that for all
$a,b,c$ and all choices of $x$:
\begin{eqnarray}
d(x,a) \leq \|x'-a'\|.
\label{eq:cat_ineq}
\end{eqnarray}
For a CAT(0) space (i) there is a unique geodesic between any two points, (ii) the space is contractible, in the topological sense, to a point and (iii) the intrinsic mean in terms of the geodesic distance is unique.
See \cite{gromov-1987} for properties of CAT(0).

Next consider CAT($k$) space which in essence generalizes CAT(0) space.
Consider a geodesic triangle $abc$ whose perimeter is less than $2\pi/\sqrt{\max(k,0)}$
for $k\in \mathbb{R}$ and a comparison triangle $a'b'c'$ on a surface $\mathcal{M}_k$ with a constant curvature $k$.
If the inequality (\ref{eq:cat_ineq}) holds for $x$ and $x'$ selected in the same manner but 
with the geodesic length $d_{\mathcal{}M_k}(x',a')$ on the surface $\mathcal{M}_k$ is used instead of 
the Euclidean distance $\|x'-a'\|$, we say the geodesic metric space has CAT($k$) property.
Thus every CAT($k$) space is a CAT($k'$) space for $k<k'$.
Intuitively speaking, CAT(0) space is a space with non-positive sectional curvatures
and CAT($k$) space is a space with sectional curvatures at most $k$.
See, for example, \cite{bridson-1999} for detailed explanation of CAT(0) and CAT($k$) spaces.

\subsection{Geodesic metrics on distributions}
Let $X$ be a $d$-dimensional Euclidean random variable absolutely continuous with respect to the Lebesgue
measure, with density $f(x)$. Let $\Gamma = \{z(t), t \in [0,1]\}$ be a parametrised integrable path between two points  $x_0=z(0),x_1=z(1)$ in $\mathbb{R}^d$, which is rectifiable with respect to the Lebesgue measure.  Let
\begin{equation*}s(t) =  \sqrt{\sum_{i=1}^d\left(\frac{\partial z_i(t)}{\partial t}\right)^2},\end{equation*}
with appropriate modification in the non-differentiable case,
be the local element of length along $\Gamma$. The weighted distance along $\Gamma$ is
\begin{equation}
d_{\Gamma}(x_0,x_1) = \int_0^1 s(t) f(z(t)) dt
\label{eq:geodesic_dist}
\end{equation}
The geodesic distance is
\begin{equation*} d(x_0,x_1) =  \inf_{\Gamma}d_{\Gamma}(x_0,x_1).\end{equation*}
Here we consider a random variable on Euclidean space but
this can be generalized for Riemannian manifolds and even for singular spaces
with a density with respect to a base measure naturally defined
by the metric.

From the geodesic distances on distributions we shall follow three main directions:
\begin{enumerate}
\item transform the geodesic metrics in various ways with parameters $\alpha, \beta$ to obtain a wide class of metrics,
\item discover (locally) CAT(0) and CAT($k$) spaces for certain ranges of the parameters,
\item apply empirical versions of the metrics based on an empirical graph whose nodes are the data points.
\end{enumerate}

There is an important distinction between global transformations applied to the
whole distance between points and local transformations applied to dilate the distance element.

\section{The $d_{\alpha}$ metric and the geodesic subgraphs}\label{sec:alpha}
The general $d_\alpha$ metric is a dilation of the original distance $d$ and  what we have referred to as a local metric.
It is obtained by transforming the density in (\ref{eq:geodesic_dist}).
Thus for $\Gamma = \{z(t), t \in [0,1]\}$ between $x_0=z(0)$ and $x_1=z(1)$,
\begin{equation*}
d_{\Gamma,\alpha}(x_0,x_1) = \int^1_0 s(t) f^\alpha(z(t)) dt 
\end{equation*}
and
\begin{equation*}d_\alpha(x_0,x_1) =  \inf_{\Gamma}d_{\Gamma,\alpha}(x_0,x_1).\end{equation*}
Here $\alpha$ is any real number.
Changing $\alpha$ essentially changes the local
curvature. Roughly speaking, when $\alpha$ is more negative
(positive), the curvature is more negative (positive).
In section \ref{sub:choosing_alpha}, we will explain how to select the value of $\alpha$ for data analysis. Values between -5 to 1 are usually selected.

In the next subsection, we look at the one-dimensional case.
Although this case is elementary, good intuition is obtained by rewriting
the standard version in terms of a geodesic metric.

\subsection{One-dimensional means and medians}\label{onedim}
Assume that $X$ is a continuous univariate random variable with probability density function $f(x)$ and cumulative distribution function (CDF) $F(x)$.
The mean $\mu = \makebox{E}[X]$ achieves $\displaystyle\min_m \makebox{E}[(X-m)^2]$. Here we are using the Euclidean distance:
$d_E(x,y) = |x-y|$.

The median is  defined by $\nu = F^{-1}(1/2)$.
On a geometric basis, we can say that $\nu$ achieves
$\displaystyle\min_m E_X[d_D(m, X)^2]$, where we use a metric that measures the amount of probability between $x$ and $z$:
\begin{equation}
d_D(x,z) = |F(x) - F(z)|. \label{cdf}
\end{equation}
Carrying out the calculations:
\begin{eqnarray*}
    \mbox{E}_X [d_D(m,X)^2] &=  \int_{-\infty}^{\infty} (F(m) - F(x))^2  f_X(x) dx\\
     &=\frac{1}{3} - F(m)(1-F(m))
\end{eqnarray*}
which achieves a minimum of $\frac{1}{12}$ at $F(m) = \frac{1}{2}$, as expected.


Another approach for the median would be to take a piecewise linear approximation to $F$ which is equivalent to
having a density $\hat{f}$ that is proportional to $\frac{1}{x_{(n+1)} - x_{(n)}}$ in the interval $[x_{(n)}, x_{(n+1)})$.
Then, the metric is
\begin{eqnarray*}
\tilde{d}_2(x,z) = \int_{\min(x,z)}^{\max(x,z)} \hat{f}(y)dy,
\end{eqnarray*}
and
$\displaystyle\min_m \sum_{i=1}^n \tilde{d}_2(x_i,m)^2$ is achieved at $x_{(\frac{n+1}{2})}$ when $n$ is odd
and at $\frac{1}{2}(x_{(\frac{n}{2})} + x_{(\frac{n+2}{2})})$, when $n$ is even.

The idea of weighting intervals should provide intuition when we extend the intervals to edges on a graph, because edges
are one-dimensional.

\subsection{The $d_{\alpha}$ metric for graphs}
\label{subsec:alpha_graphs}
There are a number of options to define an empirical version of
the $d_{\alpha}$ metric, based on data.
One such option would be to produce a smooth empirical density
$f(t)$ followed by numerical integration and optimization to compute
the geodesics.
We prefer a much simpler method based on a metric graph whose vertices are the data points.
All geodesic computation is then restricted to the graph.
We list some candidates:
(1) the complete graph, (2) the edge graph (1-skeleton) of the Delaunay simplicial complex,
(3) the Gabriel graph, (4) the $k$-NN graph, etc.  (see, for example, \cite{okabe-2000} for Delaunay complex and Gabriel graph).
The discussion below applies to the complete graph or any connected sub-graph. 

For any such graph, define a version of the $d_{\alpha}$ distance just for edges,
\begin{equation*}\tilde{d}_{\alpha,ij} = d_{ij}^{1-\alpha },\end{equation*}
where $d_{ij} $ is the Euclidean distance from $x_i$ to $x_j$.
This can be explained by making a transformation
$ ds \rightarrow \frac{ds}{d_{ij}}.$
We refer to this as {\em edge regularization}.
We then apply $\alpha$ in the usual way to obtain
$\frac{ds}{d_{ij}^{\alpha}}.$
The new ``length" of each edge $e_{ij}$ is obtained by integrating this ``density" along the edge.
In this sense, $d_{ij}$ also plays the role of density estimation.
Although we need a regularization $d_{ij}^{-1/p}$ with respect to the dimension $p$ for density estimation
(see \cite{kendall-1963}), we manage the regularization by rescaling the parameter $\alpha$.
Note that $\alpha =1$ gives the unit length and $\alpha=0$ restores the original length.

Now we consider only the set of edges $E$ of the graph $G(V,E)$ as a metric space with
the metric defined by the geodesic:
\begin{equation*}\tilde{d}_{\alpha}(x_0,x_1) = \inf_{\Gamma}\sum_{(i,j) \in \Gamma} \tilde{d}_{\alpha,ij},\end{equation*}
where the infimum is taken over all (connected) paths $\Gamma$ between $x_0$ and $x_1$.
Here we will admit $\tilde{d}_\alpha$ as an approximation of $d_\alpha$.

Note that the graph is not a complete Euclidean graph with weights equal to
the Euclidean lengths of the edges, some edges may not be in
any edge geodesics between any pair of vertices.
\begin{defn}
For an edge-weighted graph $G$ with weights $\{d_{ij}\}$ on the graph, $G^*$,
which is the union of all
the edge geodesics between all pairs of vertices, is called the geodesic sub-graph (or geodesic graph) of $G$.
\end{defn}
\noindent We will see how the geodesic sub-graphs transform as the value of $\alpha$ changes.

We make an important general position assumption that the set
of values $\{d_{ij} \mid (i,j) \in E\}$ are distinct, that is there are no ties. 
We order the values using only a single suffix for simplicity: $d_1 < d_2 < \cdots < d_M$ where $M = |E|$.
For $\alpha <1$, this induces the
$\tilde{d}_{\alpha,i} (=d_i^{1-\alpha})$ values:
\begin{equation*}\tilde{d}_{\alpha,1} < \tilde{d}_{\alpha,2} < \cdots < \tilde{d}_{\alpha,M}.\end{equation*}

Now, consider the geodesics as $\alpha \rightarrow -\infty$. Recall that a circuit in a graph is a connected path that begins
and ends in some vertex and an
elementary circuit is a circuit that visits a vertex no more than once.  Consider an edge  $(i,j) \in E$ that has the following property which we call $Q$:
it is in an elementary circuit $\mathcal C$ of the graph in which all other edges have smaller values of $d_{ij}$ namely
\begin{equation*}d_{rs} < d_{ij} \mbox{ for } (r,s) \in \mathcal C, \; (r,s) \neq (i,j).\end{equation*}
Then, the path $\Gamma(i,j)$ (within the circuit) from $x_i$ to $x_j$  not containing the edge $(i,j)$ has length
smaller than $\tilde{d}_{\alpha,ij}$ when $\alpha$ is sufficiently negative:
\begin{equation*}\sum_{(r,s) \in \Gamma(i,j)} d_{rs}^{1-\alpha} < d_{ij}^{1-\alpha}\end{equation*}

From this argument, we see that for sufficiently large $|\alpha|$ as 
$\alpha$ approaches $-\infty$,
every edge having property Q is removed from the geodesic sub-graph,
and we obtain a tree.

Let us summarize this algorithm, which applies to a general edge-weighted graph with distinct edges. We refer to this algorithm as the
{\em backwards} algorithm. It clearly gives a tree.
\begin{enumerate}
\item Let $|E|=M$ and label the edges $e_1, \ldots, e_M$ in increasing order of their weights.
\item Starting with edge $e_M$, remove $e_M$ if it is in a cycle otherwise continue to $e_{M-1}$.
\item (General step) Continue downwards at each stage removing an edge if it is in a cycle of the remaining
subgraph.
\item Stop if no more edges can be removed using step 3.
\end{enumerate}
There is a natural {\em forwards} algorithm that also yields a tree as follows.
\begin{enumerate}
\item Let $|E|=M$ and label the edges $e_1, \ldots, e_M$ in increasing order of their weights.
\item Starting with $e_1$, add an edge if adding it does not create a cycle.
\item (General step) Continue adding an edge at each step provided 
that the addition does not create a cycle.
\item Stop if no more edges can be added.
\end{enumerate}
We have the following theorem (the proof is in the appendix).
\begin{thm}
\label{thm:mst}
Given a connected edge-weighted graph $G(V,E)$ with distinct edge weights $\{d_{ij},\; (i,j) \in E\}$, the backward and forward algorithms yield the same tree, which we call $T^*(G)$. 
Furthermore, $T^*(G)$ becomes the minimum spanning tree of $G$.
\end{thm}

For sufficiently negative $\alpha$, the tree $T$ itself, that is the tree as a metric space with metric $d_{\alpha}$, is a CAT(0) space (\cite{deza-deza-2009}).
We need to extend the metric somewhat so that it applies to the edges, in addition to the nodes. Thus, for any two points
$x,x'$ on the tree, define
\begin{equation*}d_{\alpha}(x,x') = \inf_{\Gamma(x,x')} \int_{\Gamma(x,x')} w(s)ds,\end{equation*}
where the integral is taken along the (unique) path $\Gamma(x,x')$ on the tree
and $w(s) = \frac{1}{d_{ij}^{\alpha}}$ when line element $ds$ is in edge $e$ in $\Gamma(x,x')$.
Since every metric tree is a CAT(0) space, the following is an immediate consequence of Theorem \ref{thm:mst}.
\begin{cor}\label{cat1}
There is an $\alpha^*$ such that for any $\alpha \leq \alpha^*$,
the geodesic sub-graph becomes the minimal spanning tree $T^*(G)$ endowed with the
$d_{\alpha}$ metric and, therefore, becomes a CAT(0) space.
\end{cor}

We see that for sufficiently negative $\alpha$, every geodesic defined with the $d_{\alpha}$ metric
lies in the tree $T^*$. In fact, although we started with a general connected graph, any graph for which
the edges can be mapped into a Euclidean interval gives a CAT(0) tree using this construction.

Furthermore, the geodesic subgraph ``shrinks'' as $\alpha$ changes away from 1.
\begin{thm}
\label{thm:alpha-chain}
Let $G_{\alpha}$ be an edge-weighted graph with distinct weights $\{d_{ij}^{1-\alpha}\}$ and let
$G^*_{\alpha}$ be its geodesic subgraph; then for any real $\alpha$ and $\alpha'$,
\begin{equation*}|1-\alpha'| > |1-\alpha| \Rightarrow G^*_{\alpha'} \subseteq G^*_{\alpha}.\end{equation*}
Here $\subseteq$ represents the inclusion of the edge sets.
\end{thm}
Proof. This follows from the consideration of geodesics.  An edge $(i,j)$ in $G$ is not in $G^*_{\alpha}$
if it is not a geodesic. In this case, there is an
alternative path $\Gamma$ from $i$ to $j$ such that
$d_{ij}^{1-\alpha}> \sum_{(r,s) \in \Gamma} d_{rs}^{1-\alpha}$.
However, this inequality is preserved if $\alpha$ is decreased, so that $1-\alpha$ is increased. Thus an edge absent from
$G^*_{\alpha}$ is absent from $G^*_{\alpha'}$. \qed

Note that, while Theorem \ref{thm:alpha-chain} holds for any real $\alpha$ and $\alpha'$,
in application $\alpha$ is usually set at most one since otherwise the ordering of the magnitude of $d_{ij}$s becomes the inverse by taking the $1-\alpha$-th power.

\subsection{$\alpha$ and CAT($k$)}
If a space is CAT(0), then it is CAT($k$) for all $k>0$.
Let $C(X,p,r):=\{x\in X \mid  d(p,x)\leq r\}$ be a geodesic disk of
radius $r \geq 0$ centred at $p\in X$.
Define the maximum radius $D_k(X,x)$ of the disk centred at $x$ as being CAT($k$),
that is
\begin{equation*}D_k(X,x):=\sup\{r\geq 0\mid X\cap C(X,x,r) ~\mbox{is CAT($k$)}\}.\end{equation*}
If $X$ is a metric graph,
$D_k(X,x)$ is the maximum radius of the disk which is centred at $x$ and
does not include a cycle shorter than $2\pi/\sqrt{\max(k,0)}$.

Consider a rescaling of $X$ such that the shortest (longest) edge length is 1, and
denote it as $\bar{X}$ for $\alpha\leq 1$ ($\alpha>1$).
\begin{thm}
If $|\alpha'-1|>|\alpha-1|$,
\begin{equation*}D_k(\bar{G}^*_{\alpha'},x)>D_k(\bar{G}^*_{\alpha},x)~\mbox{\rm for each}~k\in \mathbb{R}.\end{equation*}
\end{thm}
\proof
Because the $\alpha$-chain is increasing for $\alpha<1$,
each cycle in $\bar{G}^*_{\alpha}$ is removed one by one as $\alpha$ decreases.
Furthermore, each cycle length increases as $\alpha$ decreases
because, by the rescaling, every edge length is greater than 1 and
it increases as $\alpha$ decreases.
This gives the decreasing property of $D_k(\bar{G}^*_{\alpha},x)$ for $\alpha\leq 1$.
We can prove the result for $\alpha> 1$ similarly.\qed

By the theorem, $\bar{G}^*_{\alpha}$ becomes ``more CAT($k$)'' for a smaller $\alpha<1$.
Because rescaling of the graph does not affect the uniqueness of the intrinsic mean,
$G^*_{\alpha}$ tends to have a unique mean for a smaller $\alpha<1$.

\subsection{Geodesic subgraphs in 2-d with different  $\alpha$}
\begin{figure}[tbp]
\begin{center}
 \begin{minipage}{4cm}
  \begin{center}
   \includegraphics[width=4cm]{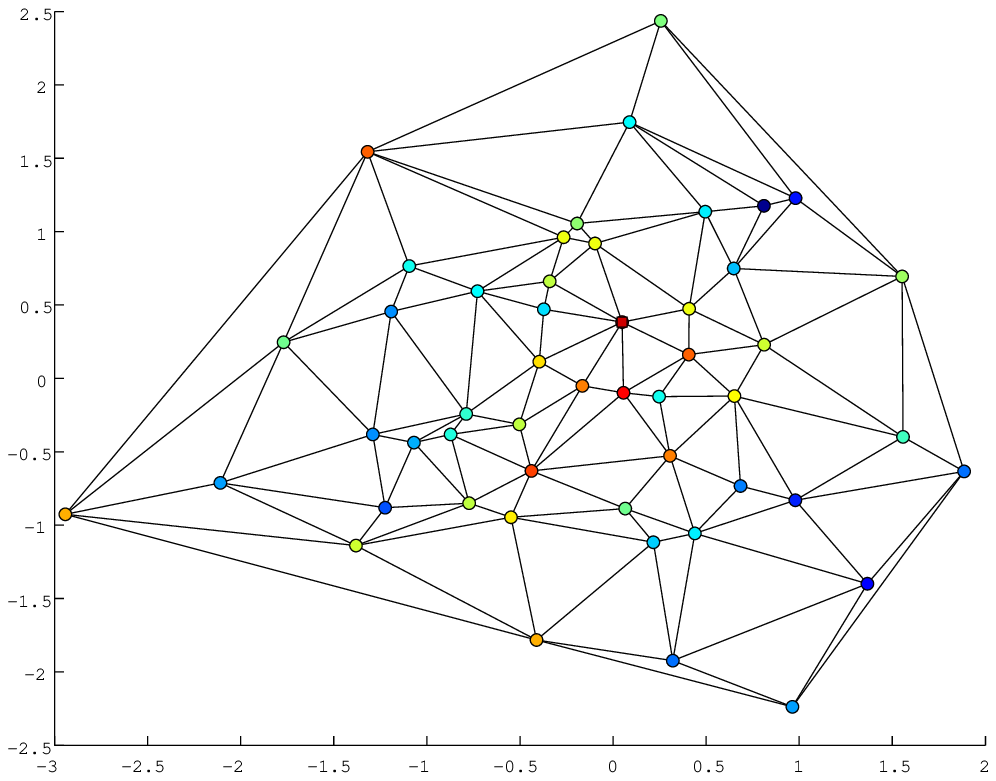}\\
   (a) $\alpha=1$
  \end{center}
 \end{minipage}
 \begin{minipage}{4cm}
  \begin{center}
   \includegraphics[width=4cm]{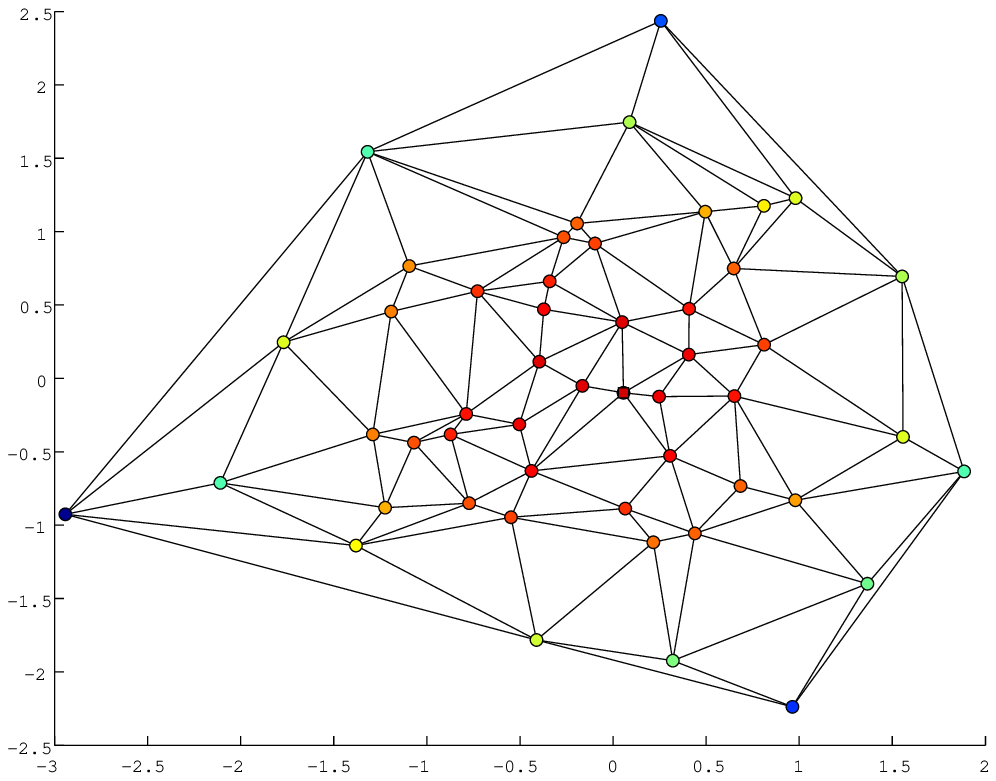}\\
   (b) $\alpha=0$
  \end{center}
 \end{minipage}
 \begin{minipage}{4cm}
  \begin{center}
   \includegraphics[width=4cm]{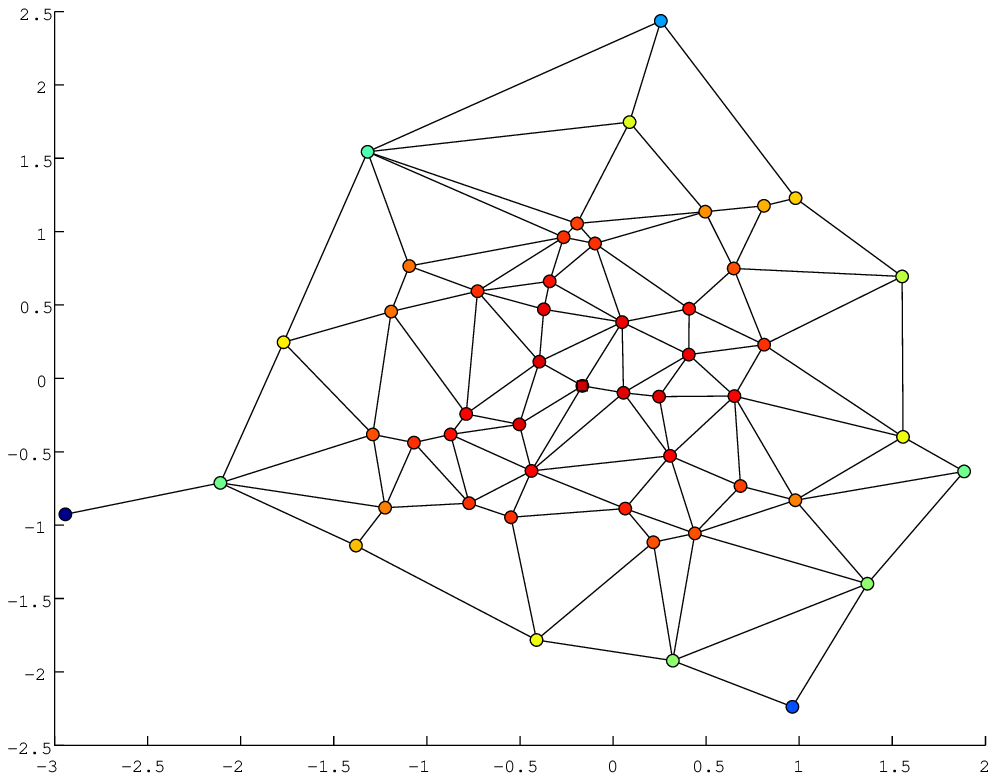}\\
   (c) $\alpha=-0.3$
  \end{center}
 \end{minipage}
 \\
 \begin{minipage}{4cm}
  \begin{center}
   \includegraphics[width=4cm]{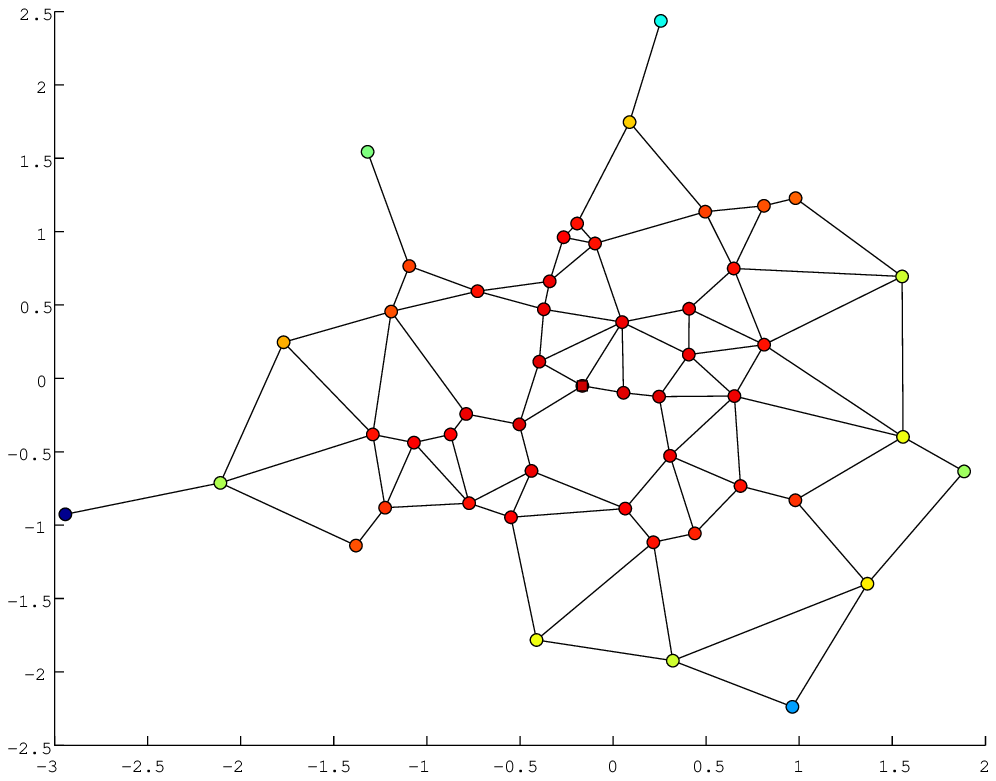}\\
   (d) $\alpha=-1$
  \end{center}
 \end{minipage}
  \begin{minipage}{4cm}
  \begin{center}
   \includegraphics[width=4cm]{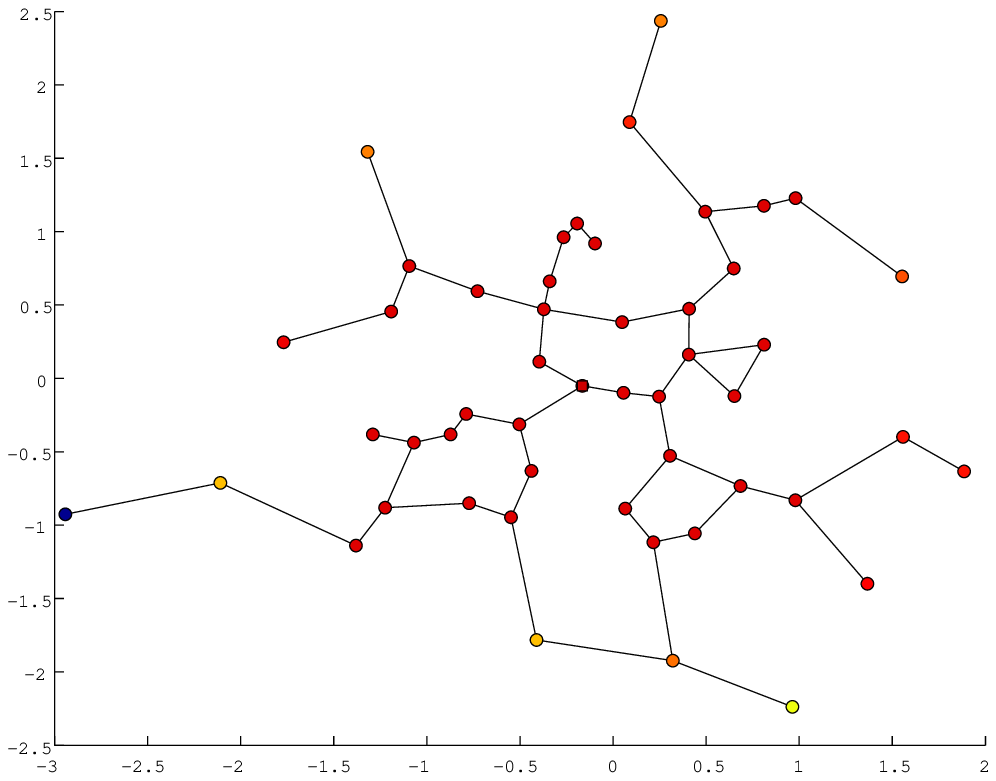}\\
   (e) $\alpha=-5$
  \end{center}
 \end{minipage}
 \begin{minipage}{4cm}
  \begin{center}
   \includegraphics[width=4.5cm]{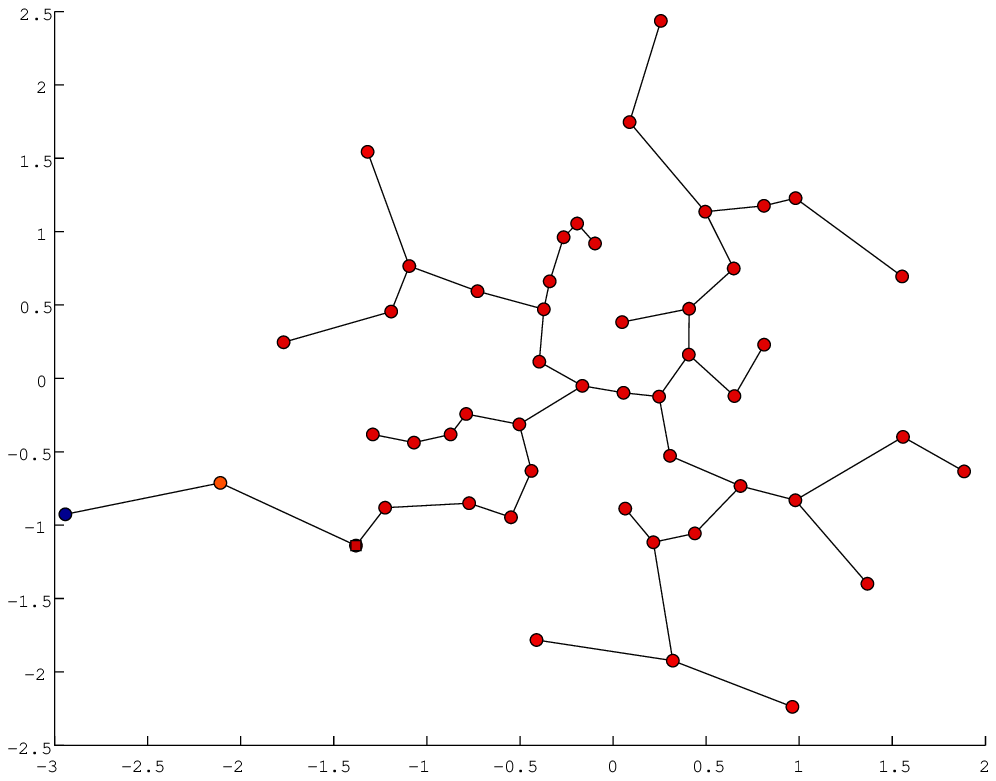}\\
   (f) $\alpha=-30$
  \end{center}
 \end{minipage}
\caption{Geodesic subgraphs with different values of $\alpha$ for 50 samples from
the standard 2-d Normal distribution. The initial graph ($\alpha=0$) is the Delaunay graph.
The value $f(x)=\sum_{i} d_\alpha(x_i,x)^2$ for each sample point $x$ is represented by the colours
red (small) and blue (large), and the minimum is represented by a square.}
\label{fig:gaussian_alpha}
\end{center}
\end{figure}

\begin{figure}[tbp]
\begin{center}
 \begin{minipage}{4cm}
  \begin{center}
   \includegraphics[width=4cm]{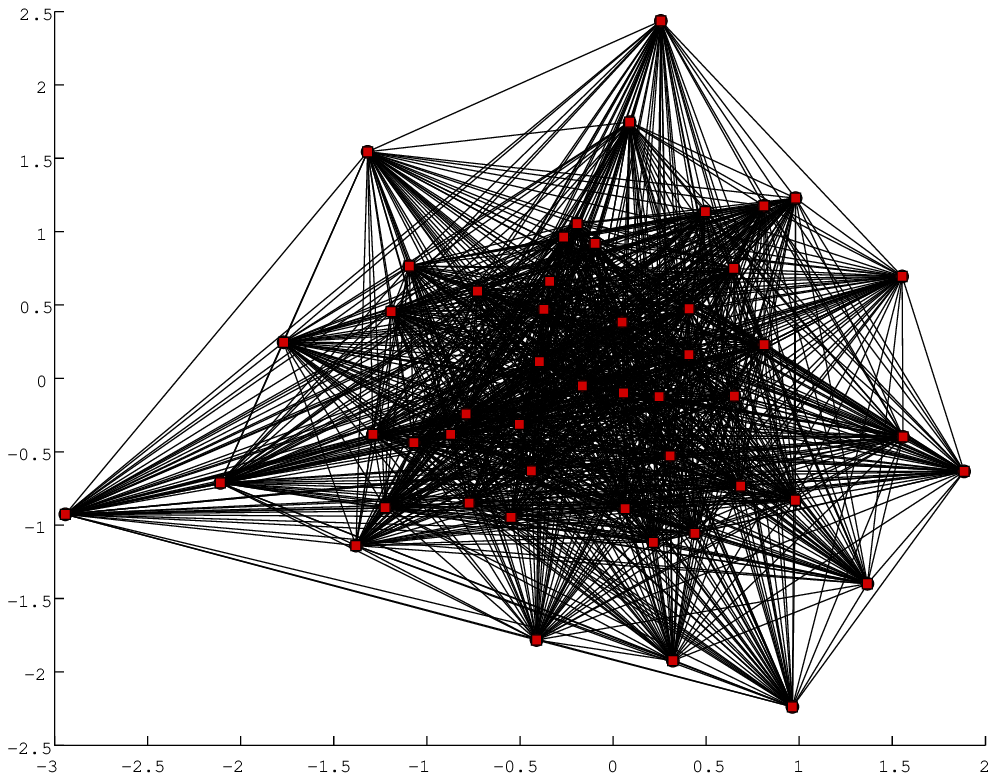}\\
   (a) $\alpha=1$
  \end{center}
 \end{minipage}
 \begin{minipage}{4cm}
  \begin{center}
   \includegraphics[width=4cm]{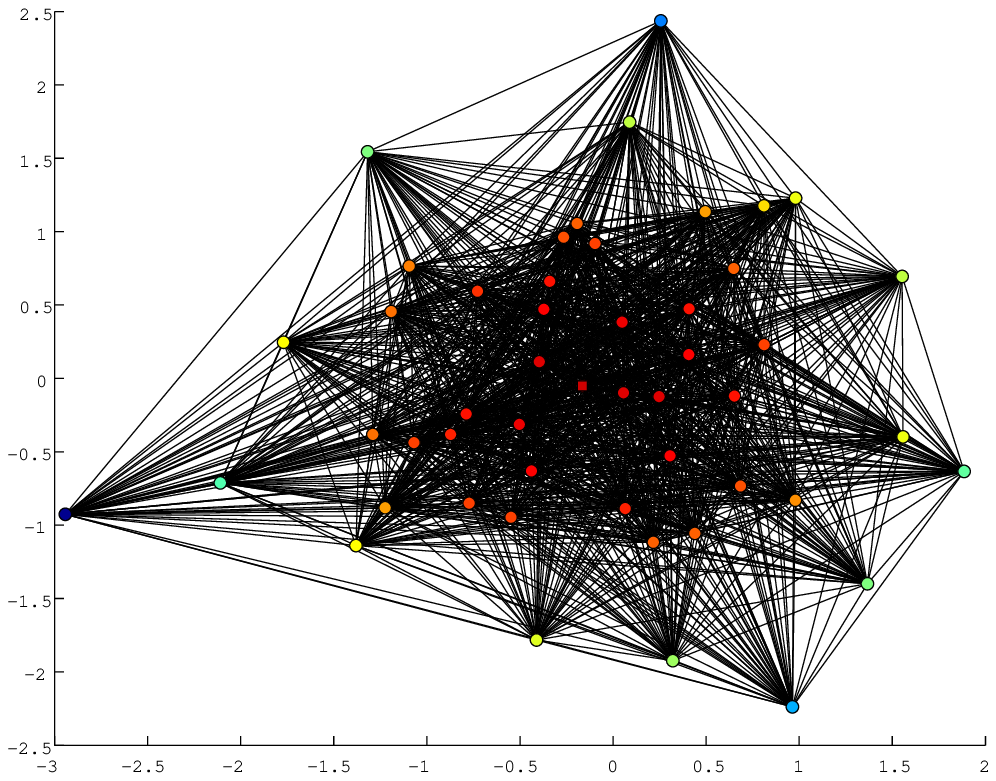}\\
   (b) $\alpha=0$
  \end{center}
 \end{minipage}
  \begin{minipage}{4cm}
  \begin{center}
   \includegraphics[width=4cm]{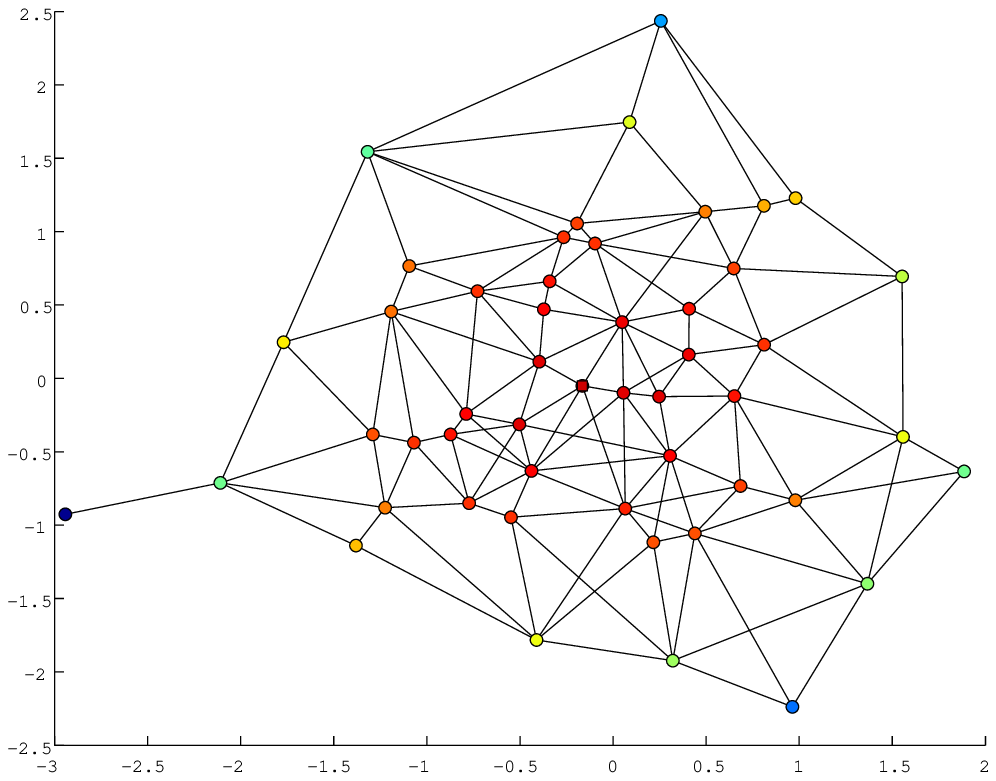}\\
   (c) $\alpha=-0.3$
  \end{center}
 \end{minipage}
 \\
 \begin{minipage}{4cm}
  \begin{center}
   \includegraphics[width=4cm]{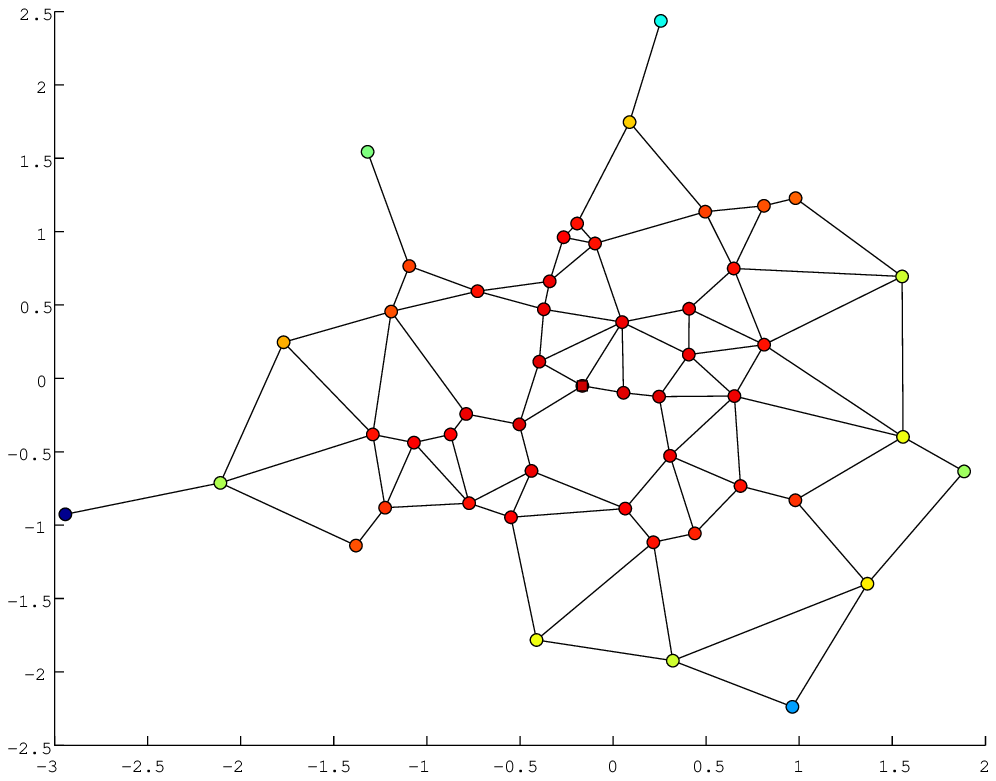}\\
   (d) $\alpha=-1$
  \end{center}
 \end{minipage}
  \begin{minipage}{4cm}
  \begin{center}
   \includegraphics[width=4cm]{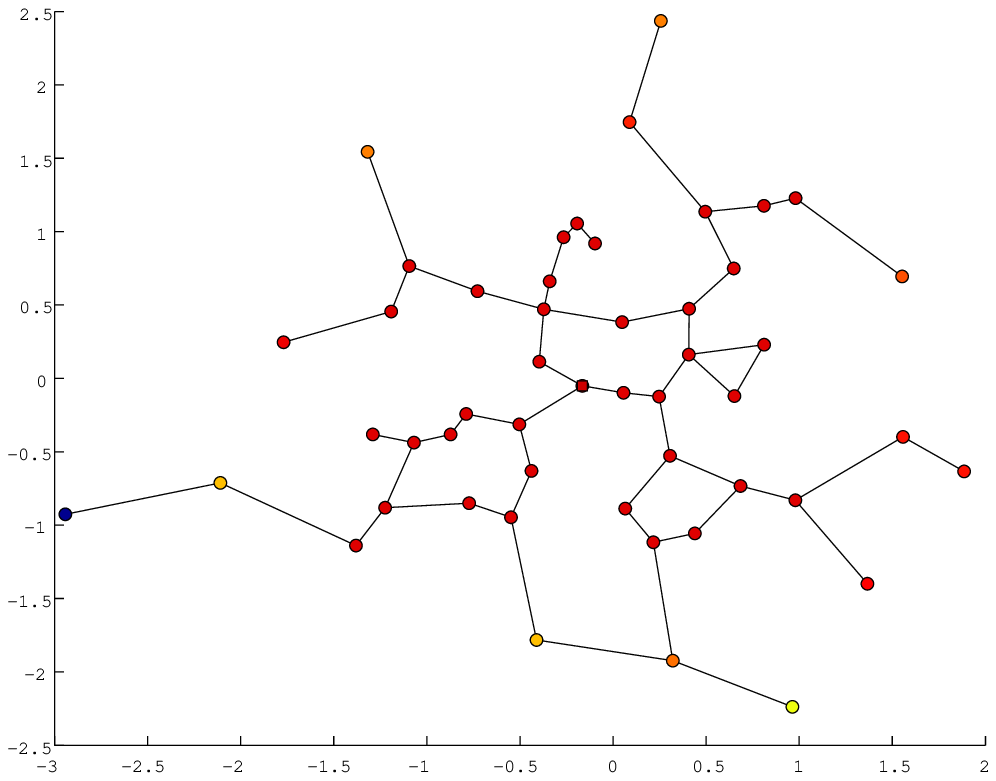}\\
   (e) $\alpha=-5$
  \end{center}
 \end{minipage}
 \begin{minipage}{4cm}
  \begin{center}
   \includegraphics[width=4cm]{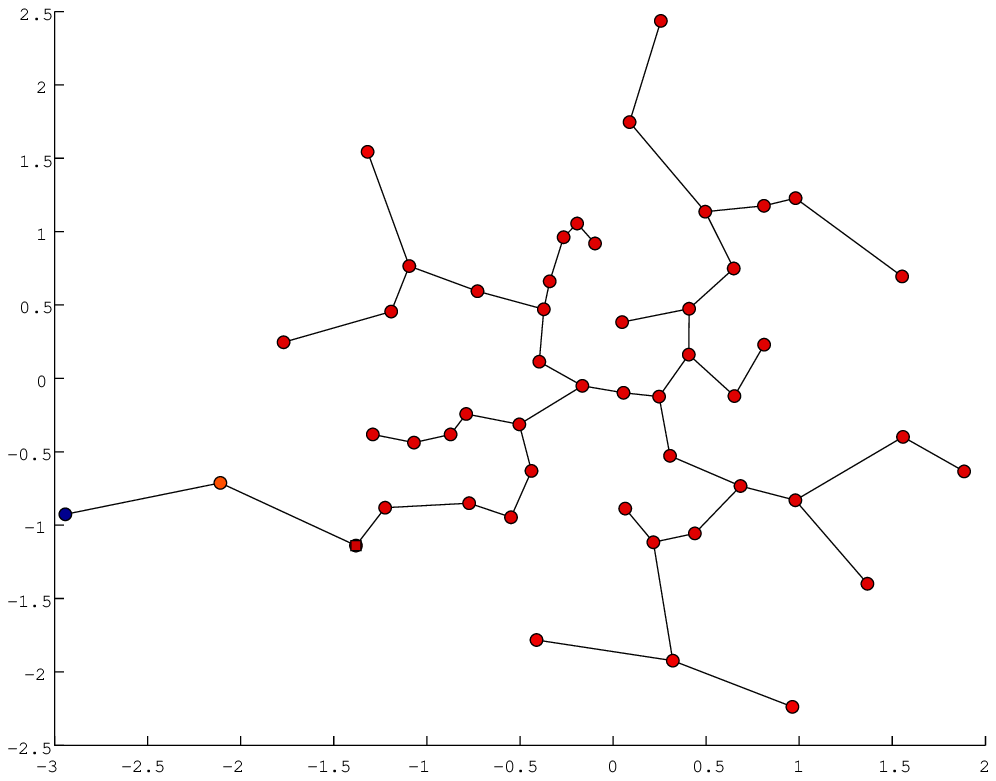}\\
   (f) $\alpha=-30$
  \end{center}
 \end{minipage}
\caption{The figures are made by the same setting of Figure \ref{fig:gaussian_alpha} but the initial graph ($\alpha=0$) is the complete graph.}
\label{fig:gaussian_alpha_complete}
\end{center}
\end{figure}

Figures \ref{fig:gaussian_alpha} (a)-(f) are geodesic subgraphs with different values of $\alpha$
for 50 samples of the standard 2-d Normal distribution. We give two cases in which we decrease $\alpha$: the Delaunay graph
in Figure \ref{fig:gaussian_alpha} and the complete graph in Figure \ref{fig:gaussian_alpha_complete}.
By the time $\alpha = -0.3$ the cases are indistinguishable and have the same
minimal spanning geodesic graph for large negative values of $\alpha$, as expected.

This is predictable from Theorem \ref{thm:alpha-chain}
and gives an important practical strategy: when the dimension is high and $\alpha$ is small,
use the complete graph rather than the Delaunay graph because the former requires computational cost
only proportional to $d$, whereas the computational cost of the latter is
$O(n^{d/2})$ (see \cite{deBerg-2008}).

\section{The $d_\beta$ metric and the metric cones}\label{sec:beta}
The CAT(0) property of Euclidean space implies that we do not obtain
multiple local minima of the Fr\'echet function $f$ even for multi-modal distributions.
However, an appropriate concave transformation of the metric
can modify the base data space
making it less CAT(0).
We introduce the $d_\beta$ metric via a  transformation $g_\beta$ as a candidate.

For any geodesic metric space $(X,d)$ with metric $d(x_0,x_1)$  and a parameter $\beta > 0$, we can define the metric
\begin{equation*}d_{\beta}(x_0,x_1) = g_{\beta}(d(x_0,x_1))\end{equation*} where
\begin{equation*}g_{\beta}(z) =  \left\{\begin{array}{l}
           \sin(\frac{\pi z}{2\beta}),\; \mbox{for} \; \; 0 \leq z \leq \beta, \\
            1, \;\; \mbox{for}\;\; z > \beta.
         \end{array}
\right.
\end{equation*}
Since $g_\beta$ is a concave function on $[0,\infty)$, 
$d_\beta$ becomes a metric but not necessarily a geodesic metric.
We can express this conveniently as $g_{\beta}(z) = \max\left(\sin(\frac{\pi z}{2\beta}),H(z-\beta)\right)$, where $H$ is the Heaviside function.

It is easiest to consider the case that $d(x_0,x_1)$ is the Euclidean distance on the real line.
As $\beta \rightarrow \infty$ for small values of $d$, the metric behaves like $\frac{\pi d}{2\beta}$, and as $\beta \rightarrow 0$, 
it behaves like $\sin(d)$ rescaled to $(0,\beta]$. For Euclidean distances greater than $\beta$, $d_{\beta}$ returns a constant distance of unity.
The metric has the effect of downsizing large distances to unity.
Because, as will soon be seen, $d_\beta$ can be recognized as a geodesic metric of a cone embedding
$X$, we refer to the mean
\begin{equation*}\hat{\mu}_\beta = \arg\inf_{\mu\in X} \sum_{i=1}^n g_\beta(d(x_i,\mu))^2\end{equation*}
as the {\it $\beta$-extrinsic mean}.

\subsection{The $\beta$-extrinsic mean: one dimension}
%
Controlling $\beta$, as will be seen below,
controls the value of $k$ when the embedding space is considered
as a CAT($k$) space.
We have an indirect link between clustering and CAT($k$) spaces.
As $\beta$ decreases while the embedding space becomes more CAT(0)
($k$ decreasing) the original space becomes less CAT(0).
This demonstrates, we believe, the importance of the CAT($k$) property in
geodesic-based clustering.

\begin{figure}[tbp]
\begin{center}
 \begin{minipage}{3.8cm}
  \begin{center}
   \includegraphics[width=3.8cm]{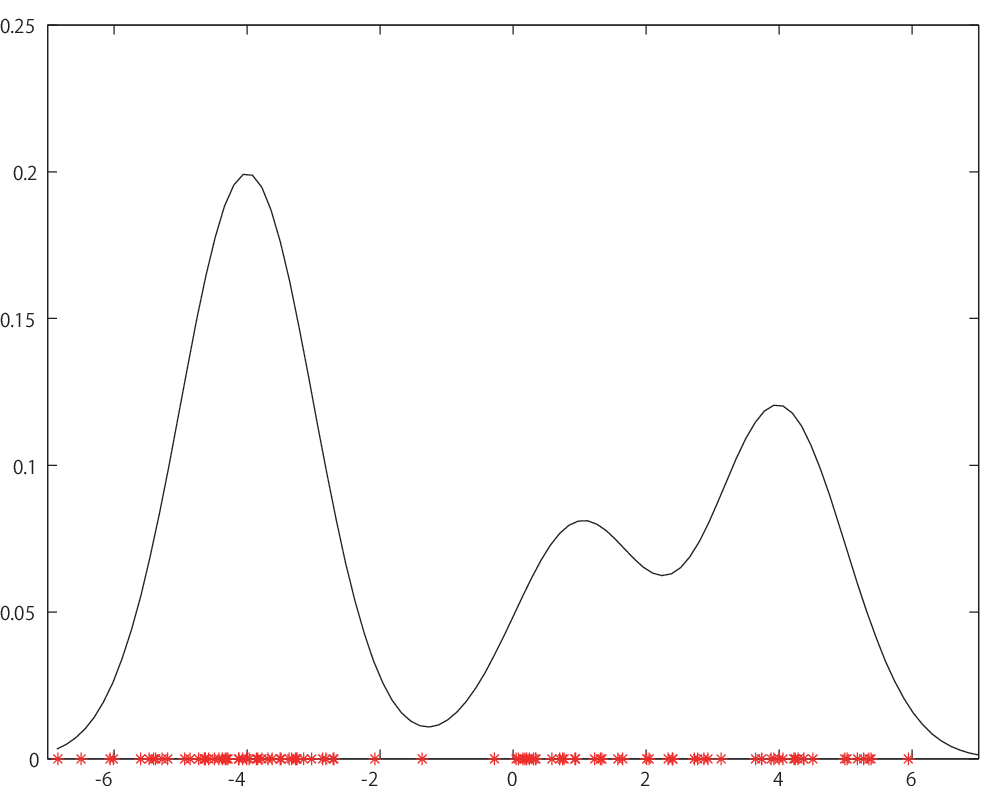}\\
   (a) Density function
  \end{center}
 \end{minipage}
 \begin{minipage}{3.8cm}
  \begin{center}
   \includegraphics[width=3.8cm]{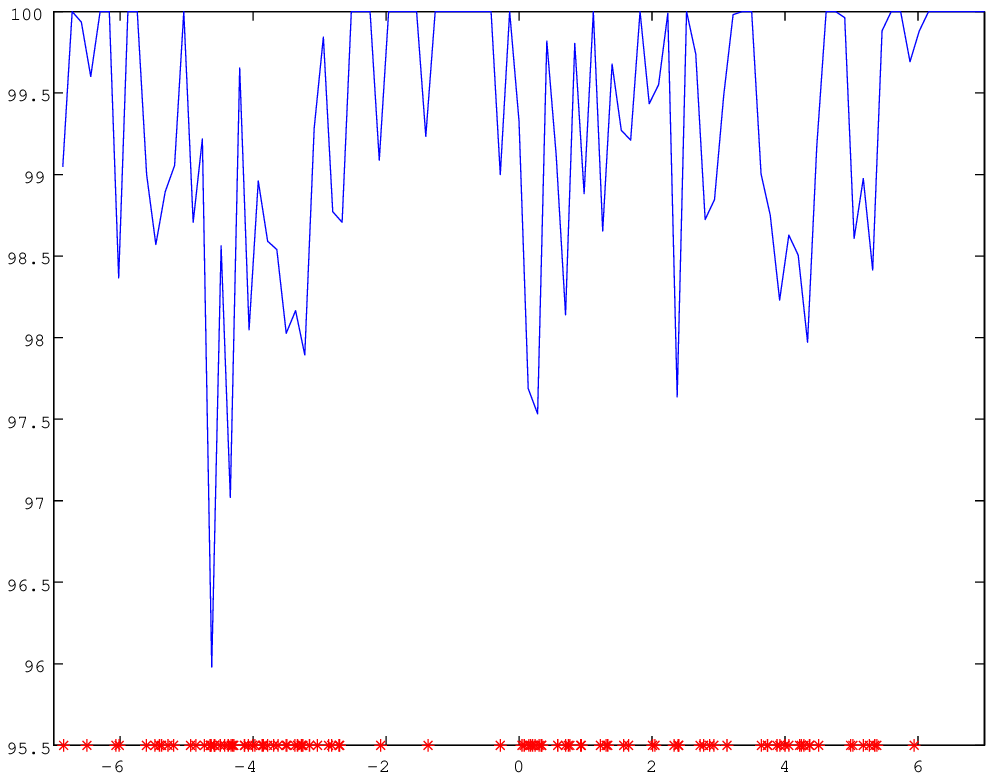}\\
   (b) $\beta=0.1$
  \end{center}
 \end{minipage}
  \begin{minipage}{3.8cm}
  \begin{center}
   \includegraphics[width=3.8cm]{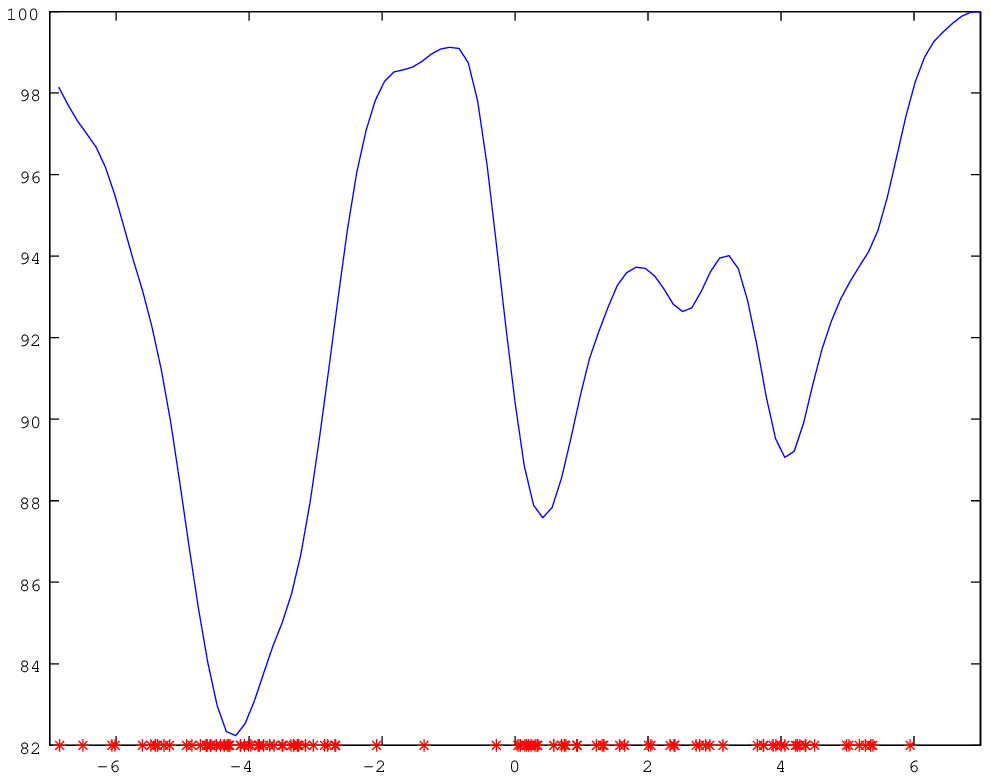}\\
   (c) $\beta=1$
  \end{center}
 \end{minipage}
 \begin{minipage}{3.8cm}
  \begin{center}
   \includegraphics[width=3.8cm]{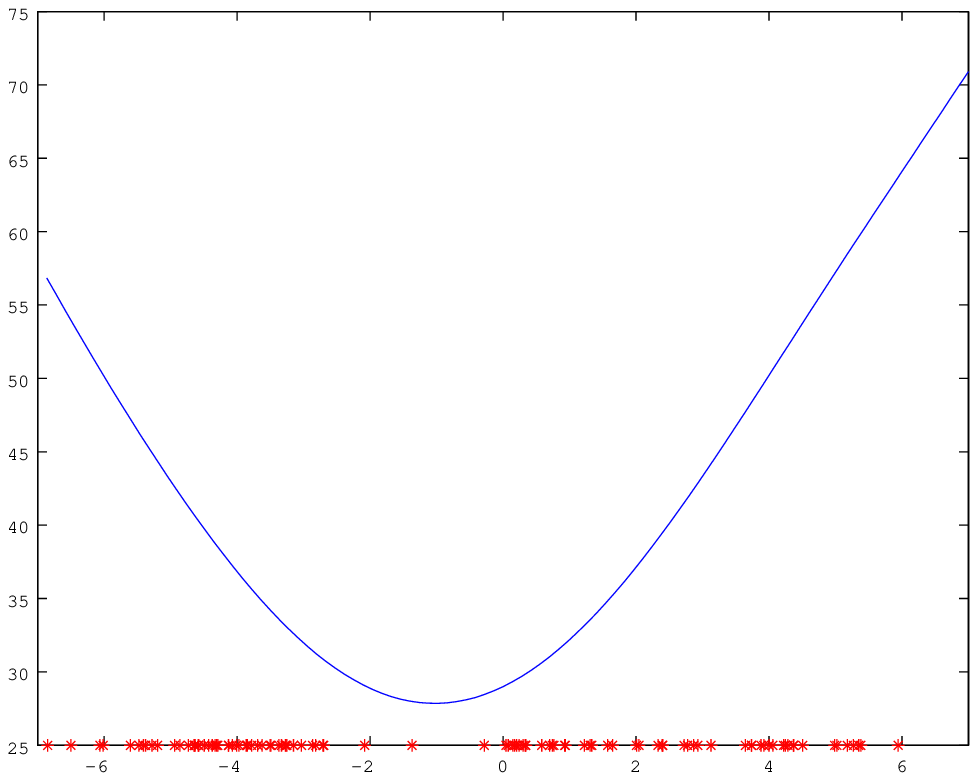}\\
   (d) $\beta=10$
  \end{center}
 \end{minipage}
\caption{(a) The density function is a mixture of three normal distributions
and 100 i.i.d. samples (red dots on the horizontal axis) from it.
(b)-(d) Graphs of $  \sum_{i} d_\beta(x_i,m)^2$ against $m$ for different values of $\beta$.}
\label{fig:1dim_beta}
\end{center}
\end{figure}

In Euclidean space, the standard Euclidean distance dose not exhibit multiple
``local means'' (i.e. local minimum points of the Fr\'echet function)
because the space is trivially CAT(0). 
However, by using the $d_{\beta}$-metric with a sufficiently small $\beta$,
the space can have multiple local means, as shown in Figure \ref{fig:1dim_beta}.

\subsection{The general case: metric cone}

The above construction is a special case of a general construction that applies to any geodesic metric space and hence to
those in this paper. Let $\mathcal{X}$ be a geodesic metric space with a metric $d_\mathcal{X}$.
A metric cone $\tilde{\mathcal{X}}_\beta$ with $\beta\in (0,\infty)$
is a cone $\mathcal{X}\times [0,1] \slash \mathcal{X} \times \{0\}$
with a metric
\begin{eqnarray*}
\lefteqn{\tilde{d}_{\beta}((x,s),(y,t))}\\
&= \frac{1}{2}\sqrt{t^2+s^2-2t s \cos(\pi\min(d_\mathcal{X}(x,y)/\beta,1))}
\end{eqnarray*}
for any $(x,s),(y,t)\in\tilde{\mathcal{X}}_\beta$.

\begin{figure}[tb]
\begin{center}
\includegraphics[height=3.5cm]{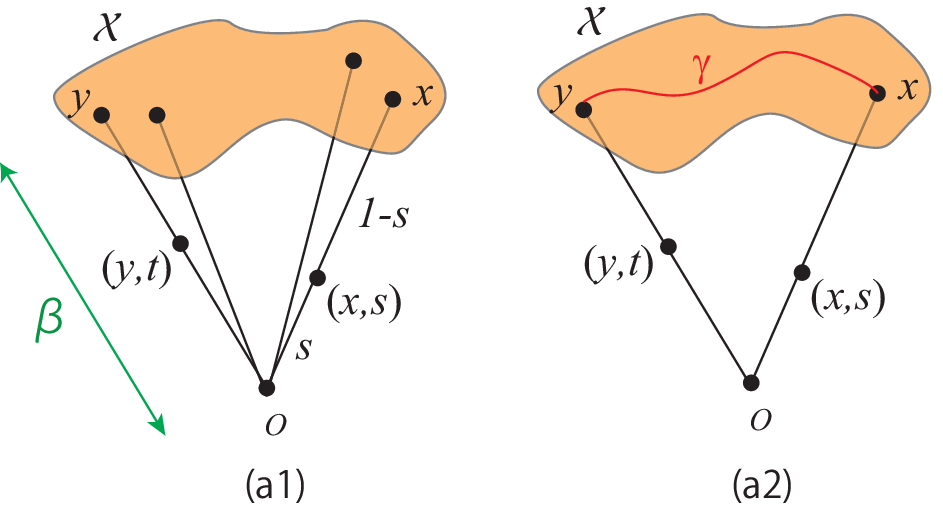}
\hspace{1cm}
\includegraphics[height=3cm]{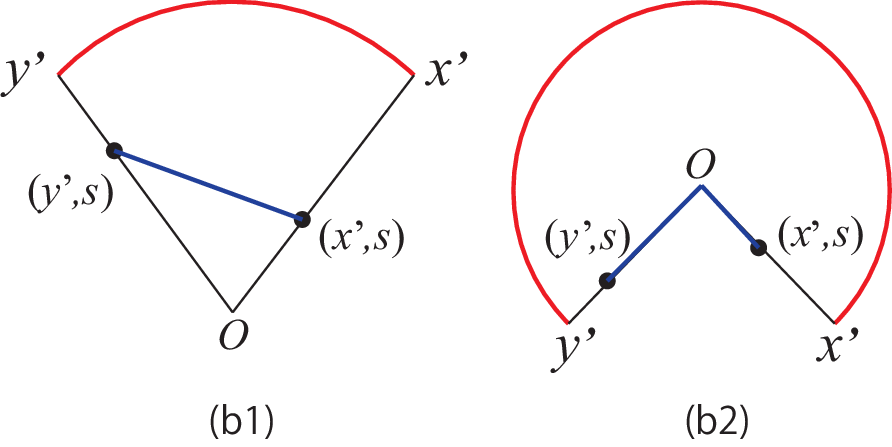}
\caption{How to define a metric cone for a geodesic metric space $\mathcal{X}$.}
\label{fig:metric-cone-explanation}
\end{center}
\end{figure}

The intuitive explanation is as follows. See Figure \ref{fig:metric-cone-explanation}.
Let $\mathcal{X}_\beta$ be the subset $\{(x,1)\mid (x,t)\in \tilde{\mathcal{X}}_\beta\}$
with the extrinsic geodesic metric on $\tilde{\mathcal{X}}_\beta$.
Thus, $\mathcal{X}_\beta$ and $\mathcal{X}$ are the same as a set but endowing different metrics.
Since $\tilde{d}_\beta((x,1),(y,1))=g_\beta(d_\mathcal{X}(x,y))$, $\mathcal{X}_\beta$ is a rescaling of the metric on $\mathcal{X}$ by $\beta$.
For any $(x,s),(y,t)\in\tilde{\mathcal{X}}_\beta$,
their projections $(x,1),(y,1)$ give two points $x,y\in \mathcal{X}$, respectively.
For a geodesic $\gamma \subset \mathcal{X}$ between $x$ and $y$, consider a
cone $\{(z,s)\mid z\in \gamma, s\in [0,1]\}$ spanned by $\gamma$.
This cone can be isometrically embedded into an ``extended unit circular sector'',
i.e. a covering $\{(r,\theta)\mid r\in [0,1],\theta \in (-\infty, \infty)\}/ \{(0,\theta)\mid \theta \in (-\infty, \infty)\}$
of the unit disk
corresponding to $\theta \in [0, \pi d_\mathcal{X}(x,y)/\beta]$.
Then $(x,s)$ and $(y,t)$ are also mapped into the extended unit circular sector;
the distance $\tilde{d}_{\beta}((x,s),(y,t))$ for $\beta=1$ corresponds
to the case (D2) of a disk if we set $(r,r')=(s,t)$ and $(\theta,\theta')=(\pi x,\pi y)$.
This corresponds to the length of the blue line path in Figure \ref{fig:metric-cone-explanation} (b1) and (b2).
For further details on metric cones, refer to \cite{deza-deza-2009}.

The following result indicates that the metric cone space preserves the
CAT(0) property of the original space and
the smaller values of $\beta$ continue this process.
\begin{thm}\label{cone}
\begin{enumerate}
\item If $\mathcal{X}$ is a CAT(0) space, the metric cone $\tilde{\mathcal{X}}_\beta$
is also CAT(0) for every $\beta\in (0,\infty)$.
\item If $\tilde{\mathcal{X}}_{\beta_2}$ is CAT(0), $\tilde{\mathcal{X}}_{\beta_1}$ is also CAT(0) for $\beta_1<\beta_2$.
\item If $\mathcal{X}$ is CAT($k$) for $k\geq 0$,
$\tilde{\mathcal{X}}_\beta$ becomes CAT(0) for  $\beta\leq \pi/\sqrt{k}$.
\end{enumerate}
\end{thm}
The proof is given in appendix \ref{proof-cone}.

It should be stressed that the theorems on $\beta$ cover metric cones based on an arbitrary geodesic metric space.
If we start with the Euclidean graph as our geodesic space, it may not be CAT(0), but it can be shown that it is a CAT($k$) space for some $k$ and will eventually  be CAT(0) for sufficiently small $\beta$.

\section{CAT($k$) spaces, curvature, diameter and uniqueness of means}\label{sec:cat_k}


In this section we prove relation between the CAT($k$) property and the uniqueness of the intrinsic means.
Let $\mathcal{X}$ be a geodesic metric space and fix it throughout this section.
The diameter of a subset $A\subset \mathcal{X}$ is
defined as the length of the longest geodesic in $A$.
We define classes $\mathcal{C}_{\rm convex}$, $\mathcal{C}_{L_\gamma}$
and $\mathcal{C}_{\rm geodesic}$ as follows.
\begin{enumerate}
\item
$\mathcal{C}_{\rm convex}$:
the class of subsets $A\subset \mathcal{X}$ such that
the geodesic distance function $f_p(x):=d(p,x)$ is strictly convex on $A$ for each $p\in A$.
Here, ``convex'' means geodesic convex, i.e. a function $f$ on $\mathcal{X}$ is convex
iff for every geodesic $\{\gamma(t)\mid t\in (t_0,t_1)\}$ on $\mathcal{X}$,
$f(\gamma(t))$ is convex with respect to $t$.

\item
$\mathcal{C}_{L_\gamma}$ for $\gamma \in [1,\infty]$:
the class of the subsets $A\subset \mathcal{X}$ such that for any
probability measure whose support is in $A$ and non-empty,
the intrinsic $L_\gamma$-mean
$$\mu=\mathop{\arg\min}_{m\in \mathcal{X}} E[d(X,m)^\gamma]$$
exists uniquely.
We refer to $\mathcal{C}_{L_2}$ as $\mathcal{C}_{\rm mean}$.

\item
$\mathcal{C}_{\rm geodesic}$:
the class of subsets $A\subset \mathcal{X}$ such that
for every pair $p,q \in A$, the geodesic between $p$ and $q$ is unique.

\end{enumerate}

\begin{lem}
\label{lem:k1-1}
\begin{equation*}\mathcal{C}_{\rm convex} \subset \mathcal{C}_{L_\gamma}
\subset \mathcal{C}_{\rm geodesic}\end{equation*}
for any $\gamma \in [1,\infty]$.
\end{lem}

\proof
If $A\in \mathcal{C}_{\rm convex}$, $f_y(x)=d(y,x)$ is a strictly convex function
on $A$ for each $y\in A$; hence, $\int d(y,x) \diff \mu_{Y}$ is strictly convex for
any probability measure $\mu$ whose support is in $A$ and non-empty.
Thus, $A\in \mathcal{C}_{L_\gamma}$.
Next, assume that $B\notin \mathcal{C}_{\rm geodesic}$ and $x,y\in B$; then, there are
at least two different geodesics, $\gamma_1$ and $\gamma_2$, between $x$ and $y$.
Thus, there are two points $x'$ and $y'$ in $\gamma_1 \cap \gamma_2$
such that there is no intersection of $\gamma_1$ and $\gamma_2$ between $x'$ and $y'$.
Then, the mid points of $x'$ and $y'$ on each geodesic become intrinsic $L_\gamma$-means
of the measure with two equal point masses on $x'$ and $y'$.
This implies that $B\notin \mathcal{C}_{L_\gamma}$.
\qed

Let $D_{\rm convex},D_{L_\gamma}$ and
$D_{\rm geodesic}$ be the largest values (including $\infty$)
such that every subset whose diameter is less than the value
belongs to $\mathcal{C}_{\rm convex},\mathcal{C}_{L_\gamma}$
and $\mathcal{C}_{\rm geodesic}$, respectively.
Then, evidently from Lemma \ref{lem:k1-1},
$
D_{\rm convex}\leq D_{L_\gamma}
\leq D_{\rm geodesic}
$
for $1\leq \gamma \leq \infty$.

Note that if $\mathcal{X}$ is CAT(0), $D_{\rm convex}=D_{L_\gamma}=D_{\rm geodesic}=\infty$.
In general, the following theorem holds.

\begin{thm}\label{diam1}
\label{thm:diameter}
\begin{itemize}
\item[(1)]
If  $\mathcal{X}$ is CAT($k$), $D_{\rm convex}\geq \pi/(2\sqrt{\max(k,0)})$.
\item[(2)]
If $\mathcal{X}$ is CAT($k$), $D_{\rm geodesic}\geq \pi/\sqrt{\max(k,0)}$.
\item[(3)]
If $\mathcal{X}$ is a surface with a constant curvature $k>0$,
$D_{L_1}\geq \pi/(2\sqrt{k})$.
\end{itemize}
\end{thm}

Some parts of Theorem \ref{diam1} are know results. See appendix \ref{proof-diameter} for details.
The proof is also given in appendix \ref{proof-diameter}.
By Theorem \ref{thm:diameter}(1), $D_{L_\gamma}\geq D_{\rm convex}\geq \pi/(2\sqrt{k})$.
Thus, a lower curvature $k$ gives a wider area where the intrinsic $L_\gamma$-mean
is unique.
According to Theorem \ref{thm:diameter}(3), this lower bound for $D_{L_1}$ is the best
universal upper bound
for any $\mathcal{X}$ with CAT($k$) property.

\vspace{0.5cm}

For $\gamma>1$, $D_{L_\gamma}$ is bounded above by
$(\theta_0(\gamma) +\pi/2)/\sqrt{k}$ where $\theta_0$
is an increasing function of $\gamma$ as shown in Figure \ref{fig:theta_0_graph}.
This bound is proved in appendix \ref{proof-gamma}.
The upper bound shows that the parameter $\gamma$ plays a role in controlling the
uniqueness of the mean, but it does not do so in Euclidean space,
where the $L_\gamma$-mean functions are always convex.
\vspace{3mm}

\begin{figure}[tb]
\begin{center}
\includegraphics[height=3.5cm]{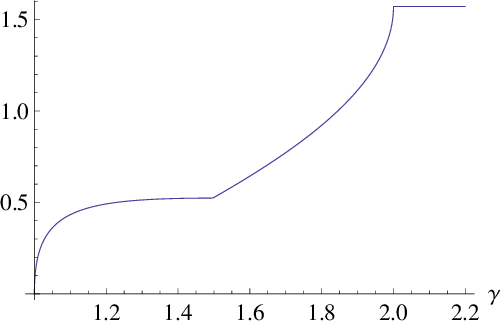}
\caption{Graph of $\theta_0(\gamma)$.}
\label{fig:theta_0_graph}
\end{center}
\end{figure}

\section{Choosing $\alpha$ and $\beta$}\label{sec:choosing}

Combining the two deformations by $\alpha$ and $\beta$, we proposed a class of
deformed metrics
\begin{equation*}d_{\alpha\beta}(x,y) = g_\beta(d_\alpha(x,y)).\end{equation*}
If we use these metrics, the Fr\'echet function becomes 
\begin{equation*}\displaystyle f_{\alpha\beta\gamma}(m)=\sum_{i=1}^n \{g_\beta(d_\alpha(x_i,m))\}^\gamma\end{equation*}
and the corresponding Fr\'echet mean and generalized variance are proposed:
\begin{equation*}\hat{\mu}=\arg\min_{m\in \M} f_{\alpha\beta\gamma}(m)~~\mbox{and}~~
{\rm Var}_{\alpha\beta\gamma}=\displaystyle\min_{m \in \mathcal{M}} \frac{1}{n} f_{\alpha\beta\gamma}(m).\end{equation*}
As explained in the previous sections, since $\alpha$ changes the curvature of the original data space and
$\beta$ changes the curvature of a metric cone embedding the data space.
Thus by tuning the values of $\alpha$ and $\beta$ we can control the uniqueness of the Fr\'echet function
via the curvatures of these two geodesic metric spaces.

In this section, we suggest how to select the values of $\alpha$ and $\beta$ empirically from the data.
For classification analysis with labels, the cross validation can be used to tune $\alpha$ and $\beta$.
Thus we will focus on the case of cluster analysis, the Fr\'echet mean and the generalized variance.

\subsection{Choosing $\alpha$} \label{sub:choosing_alpha}
First, assume that we have Euclidean data (equivalent to $\beta = \infty)$ and  recall the basic effect of decreasing $\alpha$ from $0$ to $-\infty$. At $\alpha = 0$, we  make no change to the metric.
As $\alpha$ decreases, we lose edges from the geodesic graph. 
That is to say from time to time, an edge that is in a particular geodesic is discarded and every geodesic that passes through that edge then
has to use an alternative route.

Let us assume that at $\alpha$ (and under mild extra conditions), only a single edge $e_0$ is removed and  let $d_0^{1-\alpha}$ be its length. Let $d_1^{1- \alpha}, \ldots, d_k^{1-\alpha} $ be the lengths of the edges on the  new  geodesic that will replace the removed edge.
In addition, let there be $n_0$ distinct geodesics that use $e_0$. It is straightforward to see that {\em all} geodesics that use $e_0$ will use the new arc for an interval $[\alpha, \alpha + \epsilon)$, for sufficiently small $\epsilon > 0$. The total change in geodesic length is
\begin{equation*}\Delta(\alpha) = n_0(d_0^{1-\alpha} - \sum_{i=1}^k d_i^{1-\alpha}),\end{equation*}
and it is continuous at the current $\alpha$ but the first derivative changes: $\Delta'(\alpha)$ is typically not zero. 
To see this, take the case where all the $d_i, i=1,\ldots, k,$ are equal. 
Then, the change in the first derivative is
\begin{equation*} -\frac{n_0 d_0^{1-\alpha}}{1-\alpha} \log k.\end{equation*}

In graph theory, the number of geodesics using a particular edge, $n_0$ in our case, is sometimes called the  {\em edge betweenness}. We might therefore refer to the term $n_0 d_0^{1-\alpha}$
as the weighted betweenness. This quantity measures changes in the configuration: if $n_0$ and $d_0^{1 -\alpha}$ are large then a long edge with large betweenness is removed,
and it is replaced by $k$ shorter edges from the current geodesic graph.

If $n_e$ is the betweenness of an edge $e$, the total betweenness of a graph $G(V,E)$ is the sum of all the individual edge betweennesses,
$\displaystyle\sum_{e \in E} n_e,$
and the weighted version is
$\displaystyle\sum_{e \in E} n_e d_e^{1-\alpha},$
which except for a scalar factor is  the $L_1$  variance given by $\gamma =1$, in this paper.

We shall in fact favour the use of $s_1$ ($\gamma = 2)$, and with the above discussion in mind,
we will see in Examples 1 and 2 that plots of the second derivative of $\log s_1$ do indeed have pronounced peaks and there is some matching of the $\alpha$-values at the peaks
with the analogous differential of the
aggregate betweenness.

\subsection{Choosing $\beta$}
Section 4.1 and Figure 4 are important for understanding the $\beta$ metric. We can summarise the material in a way that will indicate how to estimate $\beta$. The first point is that $\beta$ provides a metric cone. In one dimension, we wrap the real line around a circle and attach the origin. 
Then, the metric cone is based on the Euclidean metric {\em inside the cone}. 
The enlarged space (referred to as the embedding space)  is CAT(0) with respect to this metric.

We claim that this construction is fundamental because even in larger spaces,
the geodesics are one-dimensional. Every geodesic, in some sense, has its private cone but they all have a common vertex. Moreover, by Theorem \ref{cone}, if the base space is CAT(0), the embedding space is CAT(0), and in both cases, we have a unique  intrinsic mean and our statistics are well defined. However, if we compute the intrinsic mean restricted to the base space, e.g. Euclidean space, then the uniqueness no longer holds.
As stated above, the space may not be CAT(0) for small $\beta$ but may become
more so for large $\beta$. 
We can use this to our advantage: for sufficiently large $\beta$, we expect a single minimum
\begin{equation*}f_{\beta}(m) = \sum_{i=1}^n g_{\beta} (d(x_i, m))^2,  \end{equation*}
but multiple minima for smaller $\beta$, as shown in Figure \ref{fig:1dim_beta}.
If we recall that the value of the function $f_\beta$ for a given $\beta$  is helpful
in clustering, we can suggest a number of plots to show the local minima.

However, we can say more. First, note that in one dimension,
\begin{equation*}W_{\beta}(|x|)= \frac{1}{\beta}(1-g_{\beta}(|x|)^2)  = \frac{1}{2} \cos^2\left( \frac{\pi |x|}{2\beta}\right) \end{equation*}
over  $[-\beta, \beta]$ is a smooth kernel with bandwidth $\beta$.  Thus, with $d(x_i,m) = |x_i - m |$, we see that
\begin{equation*}\frac{1}{\beta} \left(1-\frac{f_{\beta}(x)}{n}     \right) = \frac{1}{n} \sum_iW_{\beta} (|x_i-x|)\end{equation*}
is a smooth density. 
This interpretation helps to intuitively choose $\beta$: select a ``typical value'' of $\|x_i-x\|$, e.g. the average of $\|x_i-x_j\|$,
by analogy with bandwidth selection for kernel functions.

Another option is to use cycle lengths in the geodesic graph. 
As can be seen from Theorem \ref{cat1} and the proof of Theorem \ref{cone}, if we set $\beta=|\Gamma|/2\pi$ for a cycle $\gamma$ and its length $|\Gamma|$, then the metric cone generated by the cycle becomes CAT(0).
The use of $\beta$ is shown in Examples 1 and 3.

\section{Examples}\label{sec:examples}
In this section, we apply the $d_{\alpha, \beta,\gamma}$ metric to real data.
Because the $L^\gamma$ loss function is more familiar than deformation of metrics by
$\alpha$ and $\beta$, we will set $\gamma=2$ and focus on
$\alpha$ and $\beta$ throughout the section.

For the $d_\alpha$ metric (for $\beta=\infty$, $\gamma=2$), we briefly describe the
computation. For each fixed $\alpha$ and each pair of points initially every
$d_\alpha(x_i,x_j)$ is computed, giving a complete graph.
On this graph the $x_i$ to $x_j$ geodesic is computed for
all $i\neq j$.
The geodesic graph, for this $\alpha$, is then computed as the union of all such
geodesics.
The present version of the software computes the geodesic graph
for a grid of around 100 points,
depending on the range of $\alpha$.
As mentioned, we are interested, here, only in the range $(-\infty, 1]$
and typically consider the range $(-r,1]$, where $r$ is a small positive integer.
For each $\alpha$, we compute $s_0^2$ and $s_1^2$.

\subsection{Example 1: $k$-nearest neighbour classification with $d_{\alpha,\beta}$}
We apply the $d_{\alpha,\beta}$ metric to the $k$-nearest neighbour ($k$-NN) method, 
one of the simplest and most popular classification methods. 
We will see that if we can choose adequate values for
$\alpha$ and $\beta$, the classification error can be reduced.

We use five data sets from he UCI Machine Learning Repository (\cite{uci-ml}):
(i) Fisher's iris data set (number of instances $n=150$, number of attributes $d=4$,
number of clusters $m=3$), 
(ii) wine data set ($n=178$, $d=13$, $m=3$), 
(iii) ionosphere data set ($n=351$, $d=32$, $m=2$), 
with only real attributes,
(iv) breast cancer Wisconsin (diagnostic) data set ($n=569$, $d=30$, $m=2$), 
and (v) yeast data set ($n=1484$, $d=8$, $m=10$). 
The average $l^2$ norm of each data set is normalized to be one.

The Euclidean complete graphs are used as the initial metric graphs ($\alpha=0$),
and classification is performed using the weighted $k$-NN method ($k$=10) 
with a common weighting
$1/d^2$ where $d$ is the distance to the neighbour data point but
using $d_{\alpha,\beta}$ for various values of
$\alpha\in \{-5,-4.8,\dots,0.8,1\}$ and $\beta\in \{2^{-6},2^{-2},\dots,2^{5},\infty\}$.
A half of the samples is selected at random as a training set
and the rest half is used as a testing set to evaluate the classification result.
We repeat it 1000 times and estimate the error rate.

\begin{table*}[htb]
\begin{center}
\caption{Classification by $k$-NN method}
\label{table:k-NN}
\begin{tabular}{|l||r|r|r|r|}
\hline
       & \multicolumn{3}{|c|}{$k$-NN with $d_{\alpha,\beta}$} & with Euclidean\\ 
\cline{2-5}
data set & $\hat{\alpha}$ & $\hat{\beta}$ & $\hat{r}$ & $r$\\ \hline

(i) iris &  -4.4  &  0.0156 &  {\bf 0.0334$\pm$0.0011}  & 0.0366$\pm$0.0011\\ \hline
(ii) wine &   0 &  $\infty$ &  0.2814$\pm$0.0025 &   0.2814$\pm$0.0025\\ \hline
(iii) ionosphere &  -0.4 &  $\infty$ &  0.1671$\pm$0.0018& 0.1677$\pm$0.0018\\ \hline
(iv)  cancer  &  0.4 &  2 &  {\bf 0.0708$\pm$0.0008}&  0.0729$\pm$0.0007\\ \hline
(v)  yeast  & 0.4  &  8 &  {\bf 0.4184$\pm$0.0009}&  0.4227$\pm$0.0008\\ \hline
\end{tabular}
\end{center}
\end{table*}

In Table \ref{table:k-NN}, $\hat{\alpha}$ and $\hat{\beta}$ are
the values attaining minimum mean classification error and $\hat{r}$ is the
error rate with 95\% confidence interval ($\pm 1.96(\mbox{std.})/\sqrt{1000}$) 
In addition, $r$ is the classification error for the ordinary Euclidean $k$-NN. 
The boldfaces represent significantly smaller error rates by $d_{\hat{\alpha},\hat{\beta}}$
than Euclidean $k$-NN.

\begin{figure*}[tbp]
\begin{center}
\begin{minipage}{5.5cm}
\begin{center}
	\includegraphics[width=5.5cm, trim={1cm 0.5cm 0 0},clip]{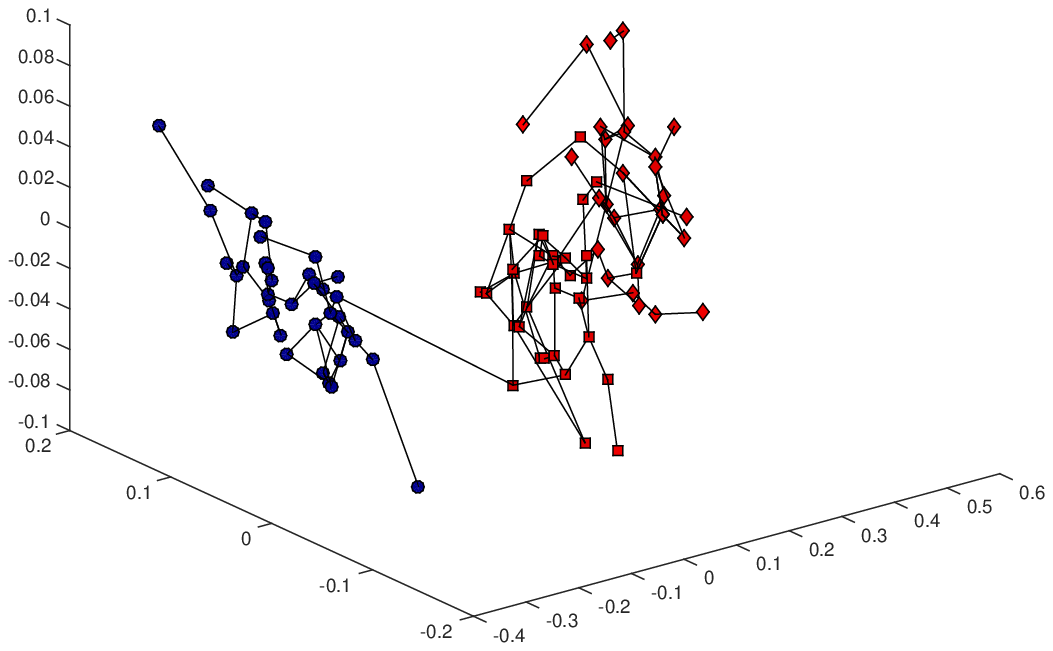}\\
   (a) iris ($\hat{\alpha}=-4.2,~ \hat{\beta}=2^{-4}$)
\end{center}
\end{minipage}
\begin{minipage}{5.5cm}
\begin{center}
   \includegraphics[width=5.5cm, trim={1cm 0.5cm 0 0},clip]{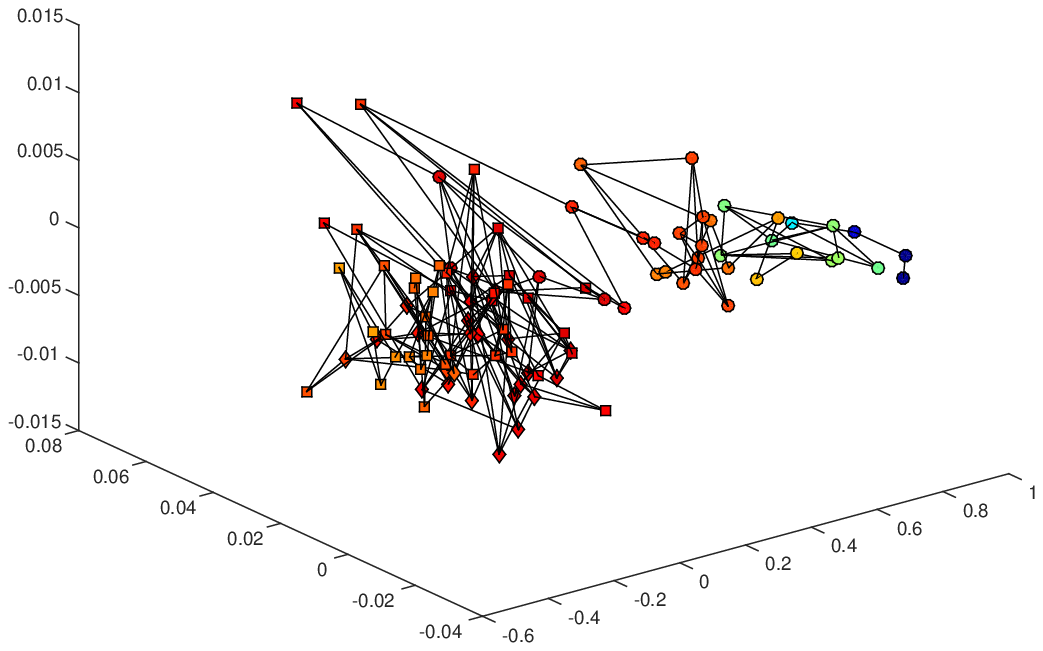}\\
   (b) wine ($\hat{\alpha}=-0.4, \hat{\beta}=\infty$)
\end{center}
\end{minipage}
\begin{minipage}{5.5cm}
\begin{center}
   \includegraphics[width=5.5cm, trim={1cm 0.5cm 0 0},clip]{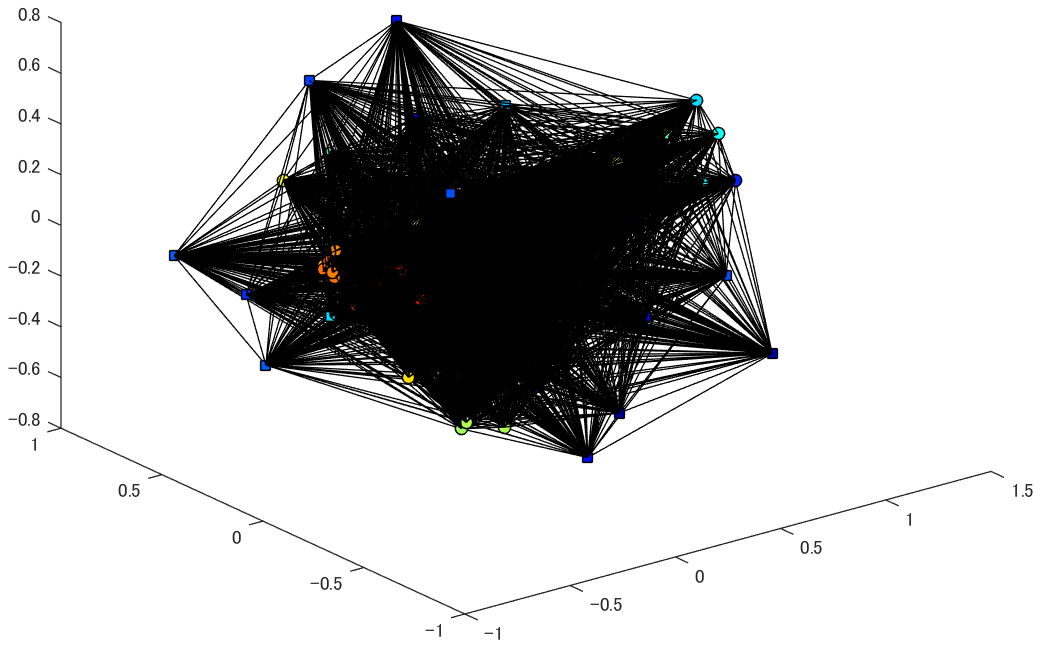}\\
   (c) ionosphere ($\hat{\alpha}=0.8, \hat{\beta}=\infty$)
\end{center}
\end{minipage}
\caption{The geodesic graph of each data set with an optimum value of $\alpha$ and $\beta$.}
\label{fig:data_geodesic_graphs}
\end{center}
\end{figure*}

Figure \ref{fig:data_geodesic_graphs} shows the geodesic graphs of the first three data sets
with the optimum values of $\alpha$ and $\beta$.
To simplify the figures, 100 samples from each data set are randomly selected
and the optimum values of $\alpha$ and $\beta$ are recomputed.
The shape of the sample points represents their class (we use only three types of point shapes 
by using the same point shape for the third and higher labeled classes for clarity of the figures).
The value of $f(x)=\sum_i g_\beta(d_\alpha(x_i,x))^2$ at each sample point is
represented by the different colours (red:small, blue:large).
We can see that the shapes of the ``optimal'' geodesic graphs are
variable because the optimal value of $\alpha$ depends on the original data spaces
and the distributions.

The computation cost is linear in the number $d$ of attributes
and therefore, the number $n$ of samples is our main concern.
The heaviest part of the algorithm is to compute the shortest path length
between each pair of samples.
We used Floyd's algorithm (\cite{floyd-1962}) which requires $O(n^3)$ computations.

There is a need for a more efficient program for more than 10,000 samples.
One option is to begin from the subgraph
of the complete graph: for example, the union of the complete subgraph
whose vertices are a subset of the samples and the edges connecting
the remaining samples to the complete subgraph.
Moreover, if we can decrease the number of edges in the
geodesic graphs, Johnson's algorithm for computing the shortest path lengths
can be used instead of Floyd's algorithm, because
it requires only $O(|E|n+n^2\log n)$, where $|E|$ is the number of edges.

\subsection{Example 2: Clustering of the world population}

\begin{figure}[tbp]
\begin{center}
 \begin{minipage}{7cm}
 \begin{center}
   \includegraphics[width=7cm]{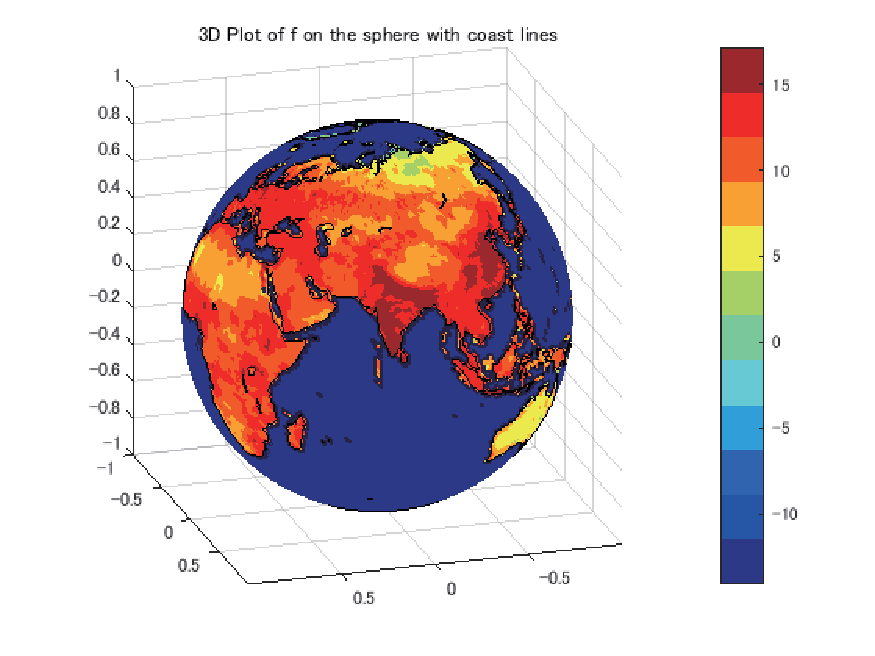}\\
  \end{center}
 \end{minipage}
 \begin{minipage}{7.5cm}
 \begin{center}
 \begin{minipage}{3.5cm}
  \begin{center}
   \includegraphics[width=3.5cm]{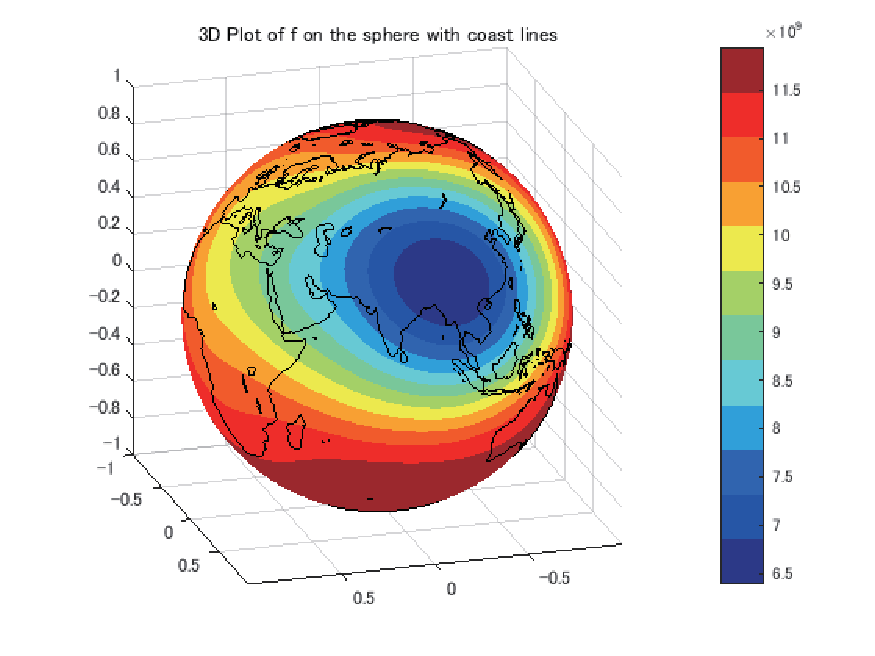}\\
   $\beta=1$
  \end{center}
 \end{minipage}
 \begin{minipage}{3.5cm}
  \begin{center}
   \includegraphics[width=3.5cm]{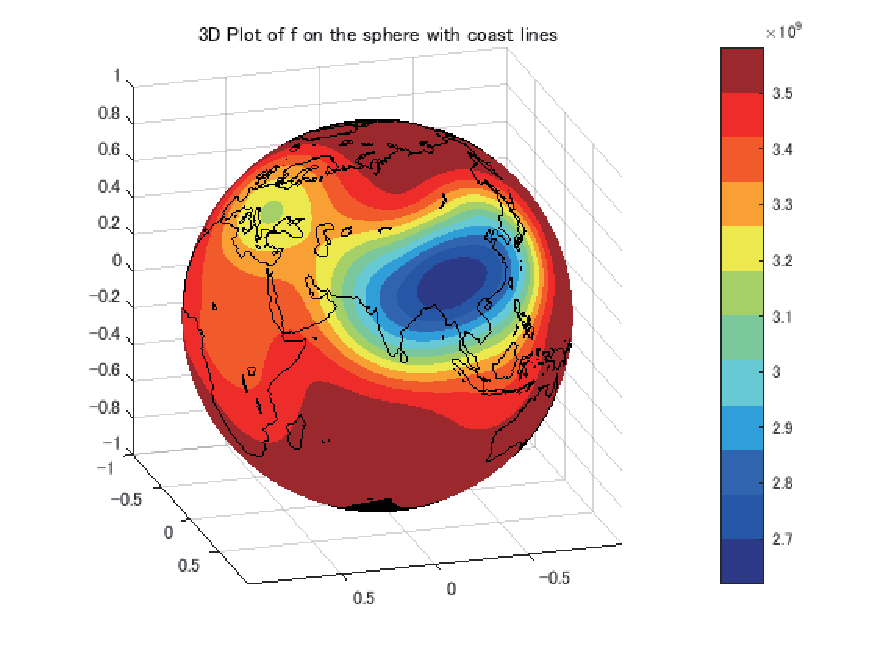}\\
   $\beta=0.3$
  \end{center}
 \end{minipage}
 \\
 \begin{minipage}{3.5cm}
  \begin{center}
   \includegraphics[width=3.5cm]{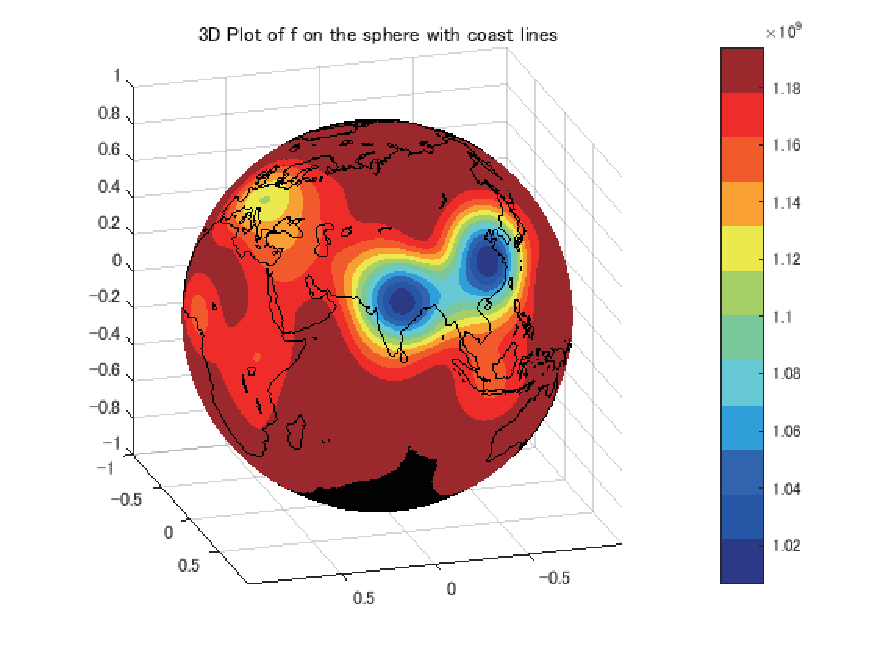}\\
   $\beta=0.1$
  \end{center}
 \end{minipage}
 \begin{minipage}{3.5cm}
  \begin{center}
   \includegraphics[width=3.5cm]{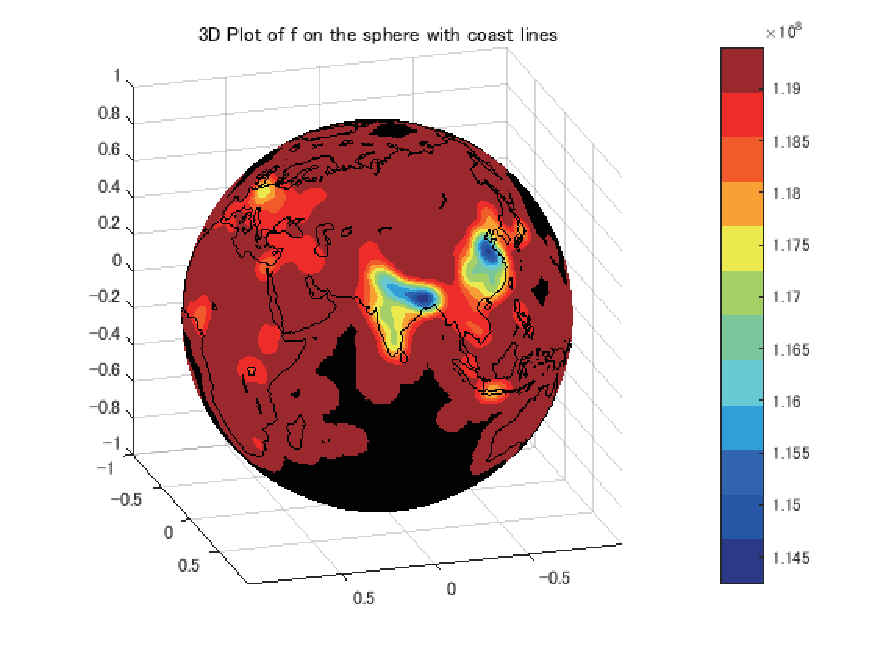}\\
  $\beta=0.01$
  \end{center}
 \end{minipage}
  \end{center}
 \end{minipage}
\caption{(left) the population density, (right) the Fr\'echet function for $\beta=1, 0.3 , 0.1, 0.01$}
\label{fig:population}
\end{center}
\end{figure}

We will show how $\beta$ plays a role in clustering analysis by using a toy example of world population.
We used the data ``Population Count Grid, v3 (2000)'' by NASA (downloadable from \cite{nasa}).
The resolution of the angle is 1 degree both for the latitude and the longitude.
Figure \ref{fig:population} (left) shows the world population density computed from the data (high:red, low:blue).
The colours in Figure \ref{fig:population} (right) represent the value of the Fr\'echet function,
\begin{equation*}f_\beta(m)=\sum_{i} g_\beta(\|x_i-m\|)^2,\end{equation*}
for $x_i,m\in S^2$.
Here the higher value of $f$ is red (lower population) and the lower value is blue (higher population).

We can see the Fr\'echet function has more local minima as $\beta$ becomes smaller.
Thus if an adequate value of $\beta$ is selected, we can obtain the centres of
a prescribed number of population clusters.
As we have seen above, a smaller value of $\beta$ corresponds to a smaller curvature
of the embedding metric cone in the sense of the CAT($k$) property.
Thus this example shows how the curvature of the embedding metric cone
affects the Fr\'echet function and the clustering analysis by the function.

\subsection{Example 3: comparison of empirical graphs via connectedness and graph Ricci curvature}

In this section, we compare the structure of empirical graphs computed by three different methods, the $\epsilon$-neighbourhood graph, the $k$-nearest neighbours graph and the $\alpha$-graph (geodesic subgraph) proposed in this paper.
The $\epsilon$-neighbourhood graph is an undirected (empirical) graph such that two vertexes are connected if the distance of the two vertexes is smaller than  a positive $\epsilon$.
The $k$-nearest neighbours ($k$-NN) graph is an undirected (empirical) graph constructed by joining each vertex to its $k$ nearest neighbour vertexes.
While the $\epsilon$-neighbourhood graph is a natural option for empirical graphs if the data points are almost uniformly distributed, the $k$-NN graph has several merits in application (e.g. the graph has usually fewer connected components) and is used more often especially for high dimensional data.
It is worth to remark that the $k$-NN algorithm is the most popular method to construct empirical graphs in the area of manifold learning.
See Section \ref{sec:survey} and the references there for more details of the manifold learning.

We use three artificial data and one real data: (1) uniform sample on $[0,1]^2$, (2) uniform sample on $S^2\subset \mathbb{R}^3$, (3) uniform sample on a subset of $H^2\subset \mathbb{R}^3$ defined by the variety $x^2+y^2-z^2=1, (-1\leq z \leq 1)$ and (4) protein data 1BUW.
Sample size for artificial data (1)-(3) is 500. (4) is a data of 3-d position of 4326 atoms in a hemoglobin protein (PDB-ID:1BUW) and downloaded from Protein Data Base(PDB) (see \cite{Berman2000-bn}).

\begin{figure}[tbp]
\begin{center}
 \begin{minipage}{4.1cm}
  \begin{center}
   \includegraphics[width=4.3cm]{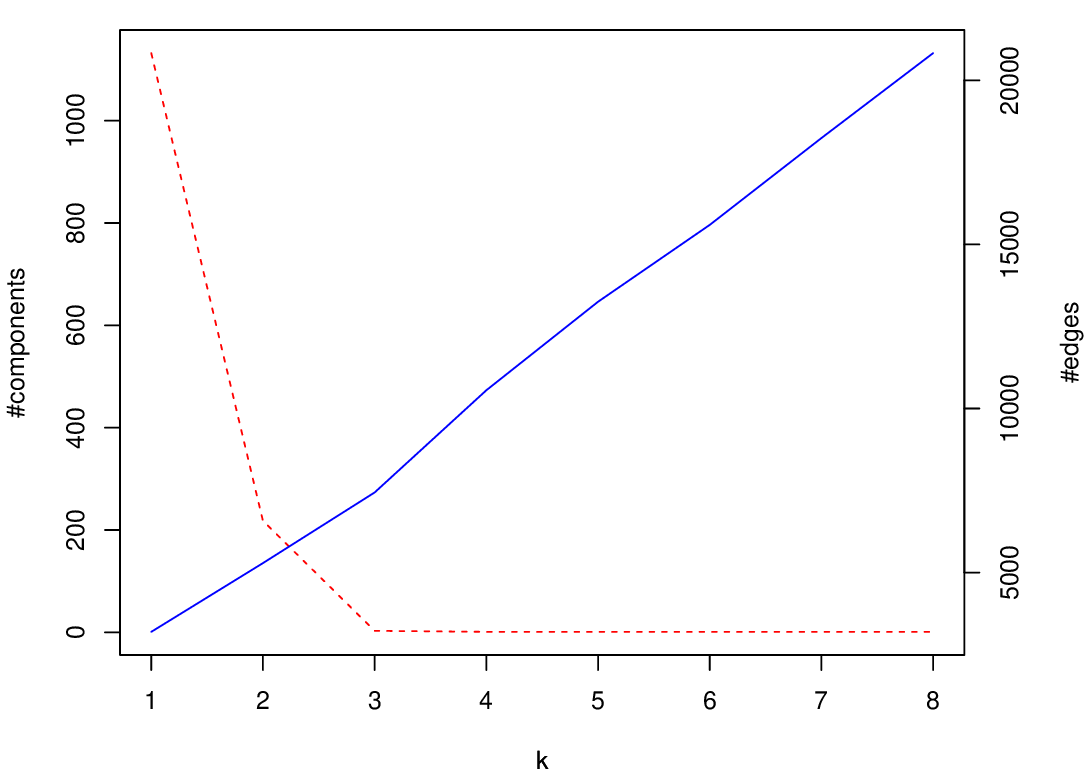}\\
   (a) $k$-NN
  \end{center}
 \end{minipage}
 \begin{minipage}{4.1cm}
  \begin{center}
   \includegraphics[width=4.3cm]{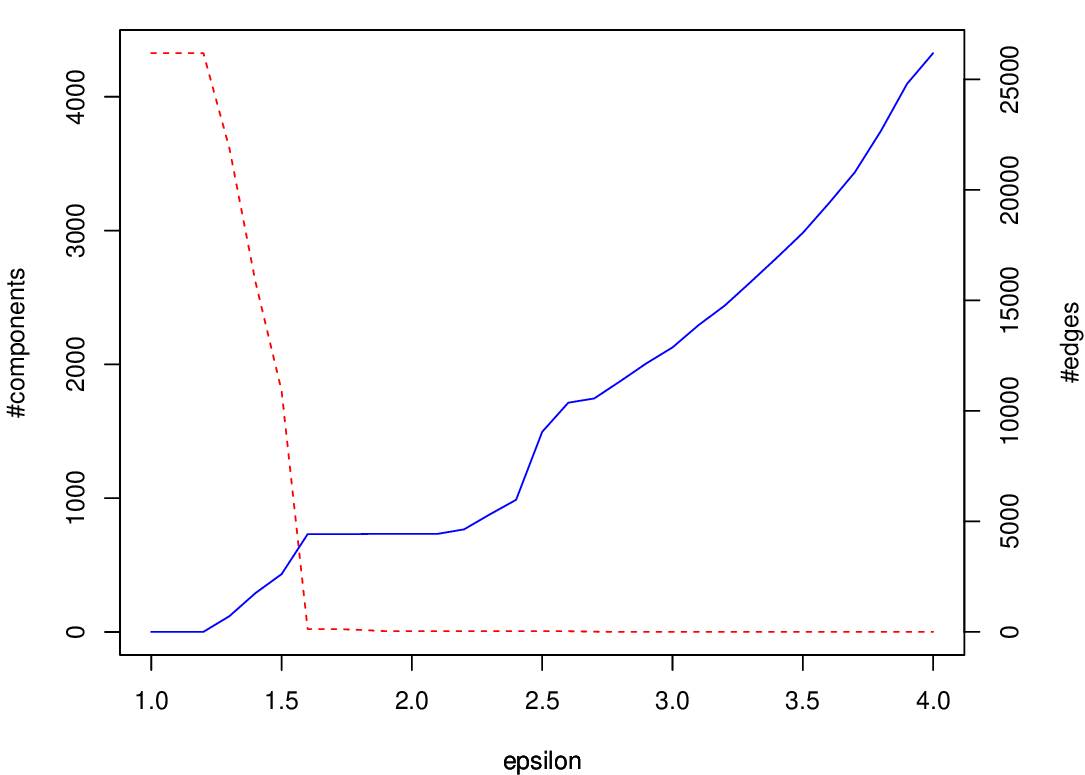}\\
   (b) $\epsilon$ neighbourhood
  \end{center}
 \end{minipage}
 \begin{minipage}{4.1cm}
  \begin{center}
   \includegraphics[width=4.3cm]{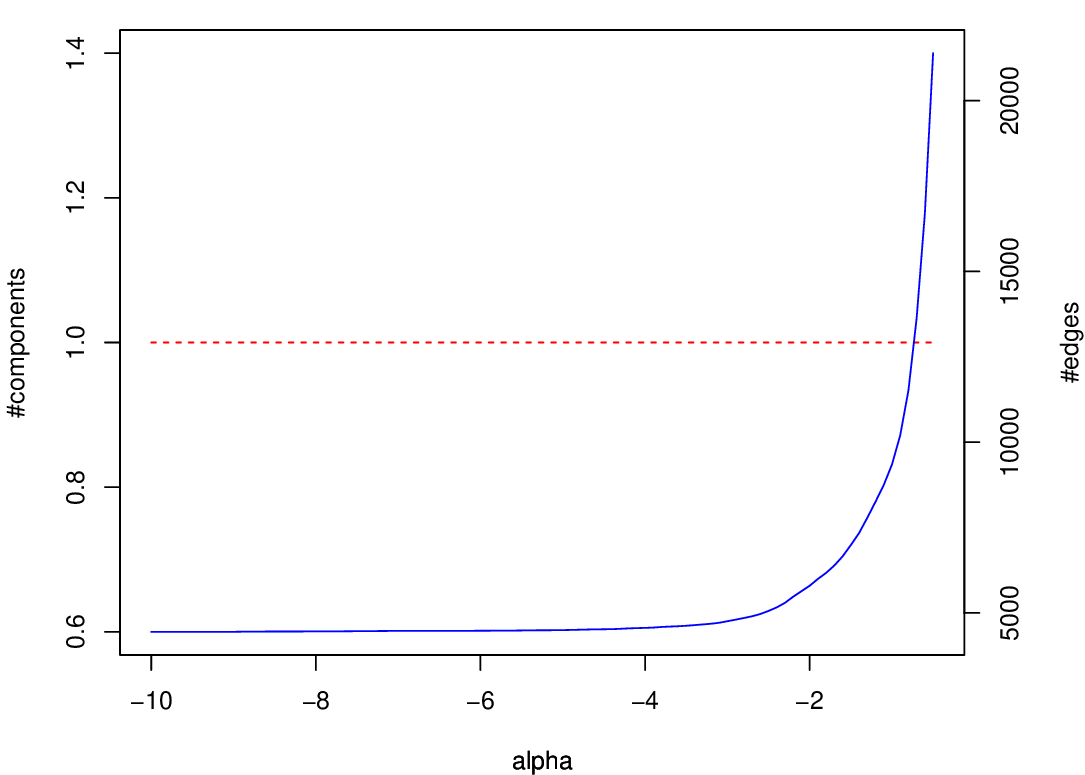}\\
   (c) $\alpha$-graph
  \end{center}
 \end{minipage}
\caption{The number of connected components (red dashed line, left axis) and edges (blue solid line, right axis) for various values of $\epsilon$, $k$ and $\alpha$, respectively, for the protein data.}
\label{fig:ex-g1}
\end{center}
\end{figure}

Table \ref{table:ex-g1} represents how the numbers of edges $e:=|E_\nu|$ and the number of connected components $c$ of the empirical graph $G_\nu:=(V, E_\nu)$ for the three graph construction methods change with parameter value $\nu$. Here, $\nu=k$ for the $k$-NN, $\epsilon$ for the $\epsilon$-neighbourhood method and $\alpha$ for the $\alpha$-graph. We first compute the $k$-NN graph for $k=1,2,4$ and $8$ and next select the values for $\epsilon$ and $\alpha$ such that the corresponding graphs have a similar number of edges.
For the $\epsilon$-neighbourhood and the $k$-NN graphs, the number of connected components changes with $\epsilon$ and $k$, respectively. As expected, the $\epsilon$-neighbourhood graph is less connected than the $k$-NN graph for non-Euclidean space.
On the contrary, the $\alpha$-graph is connected for any value of $\alpha$.
This is more evidently depicted by Figure \ref{fig:ex-g1} for the protein data.
Moreover, we can see the number of edges in the $\alpha$-graphs changes monotonically and more smoothly than other two methods when we changes the parameter value $\nu$. This means the $\alpha$ for the $\alpha$-graph is preferable for controlling the number of edges in the empirical graph.

\begin{table*}[htb]
\begin{center}
\caption{Number of edges $e$ and connected components $c$ of empirical graphs for three graph construction methods with parameter $\nu$.}
\label{table:ex-g1}
\scriptsize
  \begin{tabular}{|l|r||r|r|r|r|r|r|r|r|r|r|r|r|} \hline
    \multicolumn{2}{|l||}{}
      & \multicolumn{4}{|c|}{$k$-NN ($\nu=k$)} & \multicolumn{4}{|c|}{$\epsilon$-neighbourhood ($\nu=\epsilon$)} & \multicolumn{4}{|c|}{$\alpha$-graph ($\nu=\alpha$)} \\  \hline

%

    Uniform& $\nu$&1&2&4&8& 0.03&0.04&0.06&0.08&-10&-2.7&-0.5&-0.1 \\ \cline{2-14}
    on $[0,1]^2$& $e$& 348&652&1202&2300&333&585&1339&2366&523&649&1228&26523 \\ \cline{2-14}
    $n=500$ & $c$&152&29&1&1&247&124&6&1&1&1&1&1 \\ \hline

    Uniform& $\nu$&1&2&4&8&0.07&0.10&0.14&0.20&-10&-3.0&-0.6&-0.2 \\ \cline{2-14}
    on $S^2$& $e$& 347&648&1229&2389&340&697&1292&2470&521&643&1180&2098 \\ \cline{2-14}
    $n=500$ & $c$&153&27&1&1&318&215&118&17& 1&1&1&1 \\ \hline

    Uniform& $\nu$&1&2&4&8& 0.14&0.20&0.27&0.37&-10&-2.6&-0.5&-0.2 \\ \cline{2-14}
    on $H^2$& $e$& 350&647&1205&2370&341&666&1243&2353&519&645&1247&2322 \\ \cline{2-14}
    $n=500$ & $c$&150&24&1&1& 244&112&34&11& 1&1&1&1 \\ \hline

    Protein& $\nu$&1&2&4&8& 1.6&2.3&2.7&3.7&-10&-2.3&-0.9&-0.5 \\ \cline{2-14}
    1BUW& $e$& 3194&5289&10556&20830&4420&5321&10563&20800&4446&5295&10198&21393 \\ \cline{2-14}
    $n=4326$& $c$&1132&218&1&1&22&6&2&1& 1&1&1&1 \\ \hline
 \end{tabular}
\end{center}
\end{table*}

\begin{figure}[tbp]
\begin{center}
\includegraphics[width=8cm]{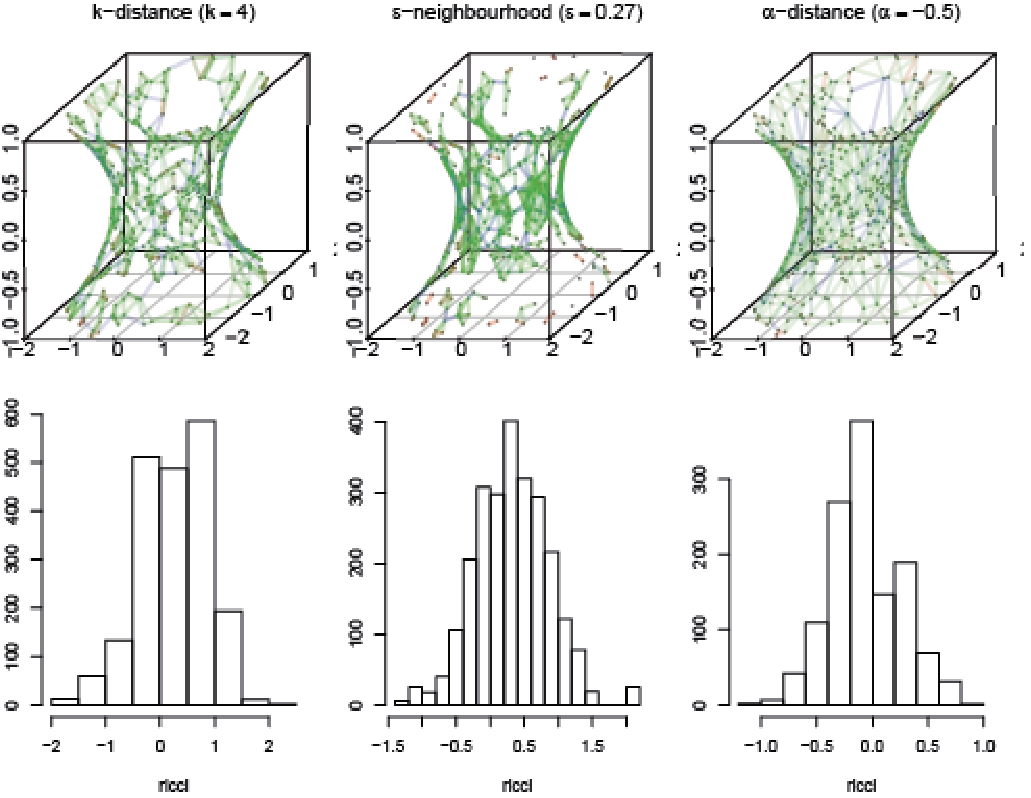}\\
\caption{The Ricci curvature of each edge (red:large, blue:small) for Uniform sample on $H^2$}
\label{fig:ex-g2}
\end{center}
\end{figure}

\begin{figure}[tbp]
\begin{center}
\includegraphics[width=8cm]{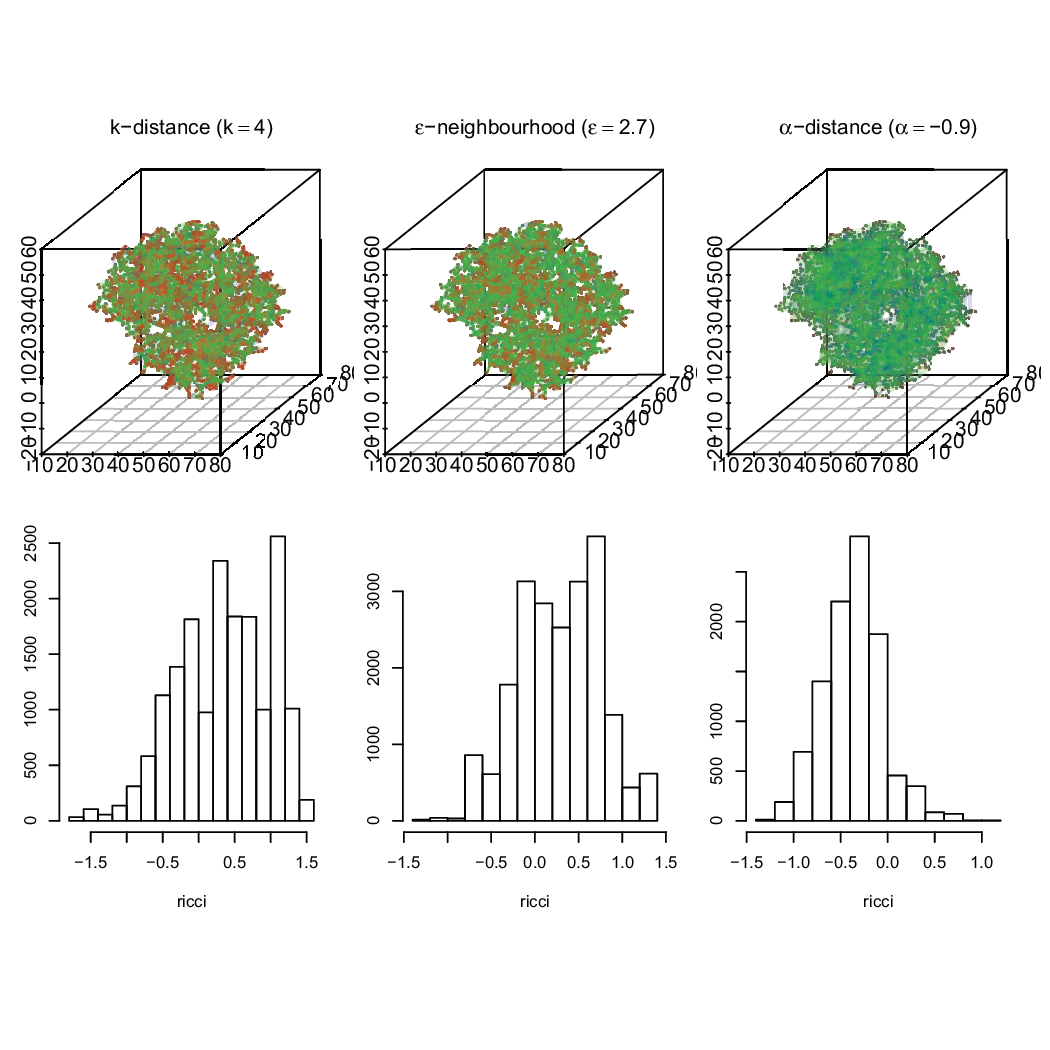}\\
\caption{The Ricci curvature of each edge (red:large, blue:small) for the protein data (1BUW)}
\label{fig:ex-g3}
\end{center}
\end{figure}

Next we compare the three types of empirical graphs via the graph Ricci curvature proposed in \cite{Lin2011-zg}. 
The graph Ricci curvature for each pair of vertexes is defined by using the graph Wasserstein metric on the graph and has some analogy to the Ricci curvature on Riemannian manifolds.
We compute the Ricci curvature for every edge in the empirical graphs for data (3) and (4).
In Fig. \ref{fig:ex-g2} for data (3) and Fig. \ref{fig:ex-g3} for data (4), each edge is coloured (blue:small, red:large) by its Ricci curvature 
for the three types of empirical graphs. Here the parameters are selected from Table \ref{table:ex-g1} as $k=4$, $\epsilon=0.27$, $\alpha=-0.5$ for Fig. \ref{fig:ex-g2} and $k=4$, $\epsilon=2.7$, $\alpha=-0.9$ for Fig. \ref{fig:ex-g3}.
The histogram of Ricci curvatures for all edges of each empirical graph is displayed under the graph.
Each histogram seems to converge to normal distribution (this is surprising for us) and the histogram for $\alpha$-graph converges faster than other two.
We expect the reason for this property is partly because the $\alpha$ controls the CAT($k$) property, another kind of curvature but related to Ricci curvature, of $\alpha$-graphs.
We remark that Ricci curvatures in the $\alpha$-graphs for $\alpha<0$ tend to have some negative bias for our examples. This is reasonable when we remember a negative value of $\alpha$ makes the data space more CAT($k$).

\subsection{Example 4: Rainfall data}
\begin{figure*}[htb]
\begin{center}
\includegraphics[height=6cm,width=12cm]{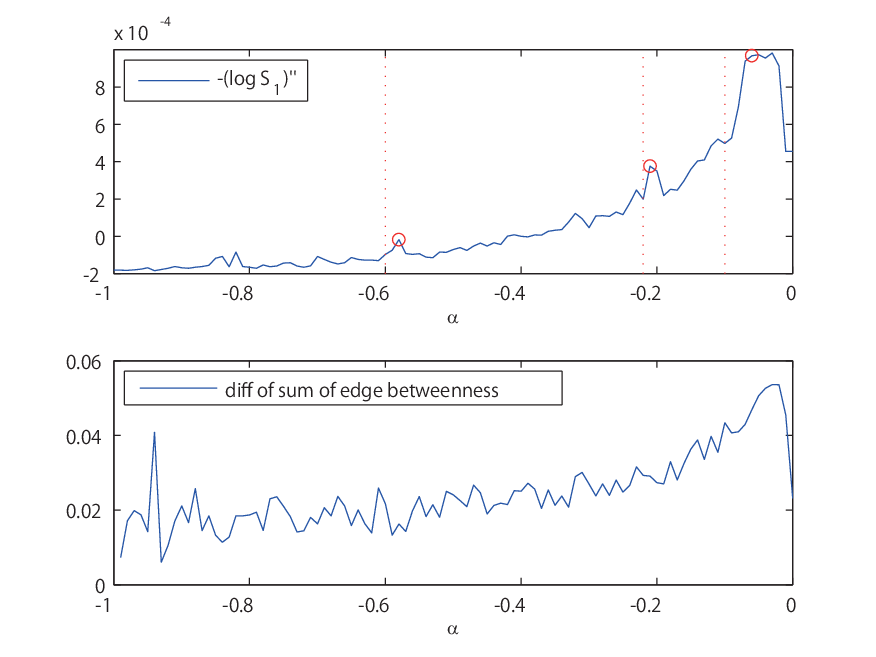}
\caption{Selecting $\alpha$ for weather data of 2014. Upper graphs: selected peaks of $-\log  {s_1}''$ (marked as circles) and values of $\alpha$ (dashed lines). Lower graphs: numerical differentiation of the sum of betweenness of each edge with respect to $\alpha$.}
\label{fig:alphaweather}
\end{center}
\end{figure*}

\begin{figure*}[htb]
\begin{center}
\includegraphics[height=5cm,width=10cm]{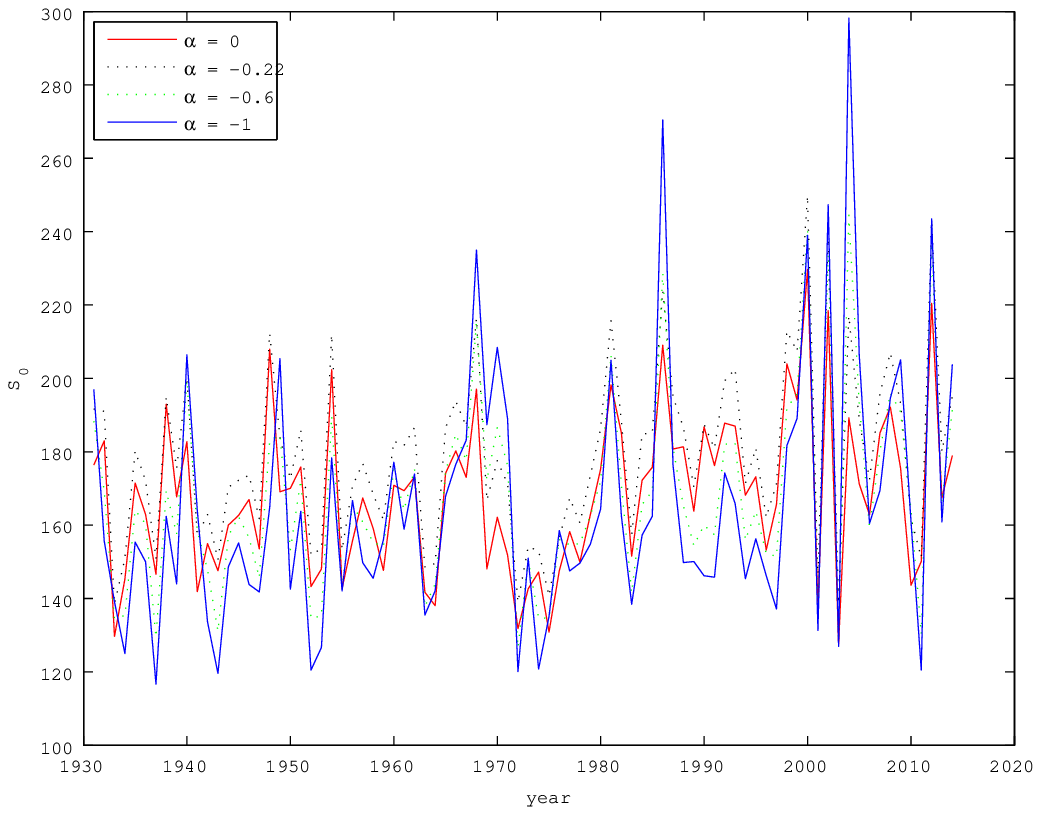}
\caption{Rainfall data in UK (1931-2014): yearly time series of $s_0^{2/(\alpha-1)}$ plotted
for $\alpha =0$ (red solid line), $-0.1$ (green dashed line), $-0.6$ (black dashed line) and
$-1$ (blue solid line).}
\label{fig:UK_weather}
\end{center}
\end{figure*}

We carry out some analysis of rainfall (precipitation) data obtained from
the UK Met Office Hadley Centre (downloadable from \cite{alexander-2001}).
Considering a single year's data we take the ``dimensions'' as the nine
regions of the UK:
South East England, South West England and Wales, Central England, North West England and Wales,
North East England, South Scotland, North Scotland, East Scotland and Northern Ireland,
and the ``points'' as the 365 (or 366) days of the year. We take the years 1931 to  2014.
Initially, we select values of $\alpha=0, -0.1, -0.22, -1$ by using some peaks
of $-(\log s_1)''$ in Figure \ref{fig:alphaweather}.
For each year, we compute $s_0^{2/(1-\alpha)}$ for values $\alpha=0, -0.1, -0.22, -1$.
The data is presented as four time series with different values of $\alpha$ for yearly values from 1931 to 2014
in Figure \ref{fig:UK_weather}. 
The figure is consistent with an emerging consensus of increased extremes and volatility in
precipitation in the UK in recent years (see Met Office, 2014,``Recent Storms Briefin'').

We now discuss the choice of $\alpha$. Following the discussion in Figure \ref{fig:alphaweather},
we plot  $- \frac{\partial^2}{\partial \alpha^2} \log s_1$ in the range $\alpha \in [-1,0]$. 
Figure \ref{fig:alphaweather}  shows plots for the year 2014. We select the values of $\alpha$ 
to be slightly smaller than the peaks. 
The local peak at approximately $\alpha = -0.22 $ indicates a rapid change in the topology of the graph at this point: we lose a considerable number of longer edges and reveal the structure in the data as a consequence. 
The betweenness plot is not so revealing, except near $\alpha = -1$. 

\section{Conclusion}\label{sec:conclusion}

The $\alpha$ metric is a deformation of the starting geodesic, and as $\alpha$ becomes
more negative, the geodesic graph, namely the union of all the geodesics, becomes sparser,
and in our graph representation, it becomes a tree.
The space is CAT($k$) with smaller $k$ and finally becomes a tree, at which point the space becomes CAT(0).
It is quite difficult to see the tree computation because of the numerous short edges, 
but for moderate values of $\alpha$ such as $-1$, 
the structure is tree-like. Abrupt changes in various statistics as $\alpha$ changes
reveal topological changes in the structure of the geodesic graph, a fact that can be used to tune $\alpha$.

The $\beta$ metric is ``non-geodesic" because although the function $g_{\beta}$ operates on a geodesic, that does not
mean that the space is a geodesic space in the formal sense. 
However the cone construction yields a geodesic metric space, which is CAT($k$)
with a lower value of $k$ than the original space, and indeed may be CAT(0). 
If the $\beta$ metric is projected back to the original space,
that space can have a non-convex Fr\'echet function with larger $k$. 
This is useful for finding clusters because of multiple minima of the Fr\'echet function, which is itself similar to a kernel. 
The means  obtained by the $\beta$-metric may represent the first study of an extrinsic mean via embedding in non-Euclidean
spaces and the first application of metric cones to statistics and data analysis.

We believe that the curvature of the  data space underlying this work demands further investigation
whereby connections should be established with
recent developments related to empirical geodesic graphs, for example in manifold learning. 
One important direction should be the effect of the curvature of the space on
the trade-off between the uniqueness of the Fr\'echet means
and the robustness of estimation. 
To this end, $\alpha, \beta$ and $\gamma$ can be seen as parameters that can be tuned to change
the curvature and hence study the trade-off.


\appendix

\section{Proof of Theorem \ref{cone}} \label{proof-cone}
\begin{figure}[tb]
\begin{center}
\includegraphics[height=4cm]{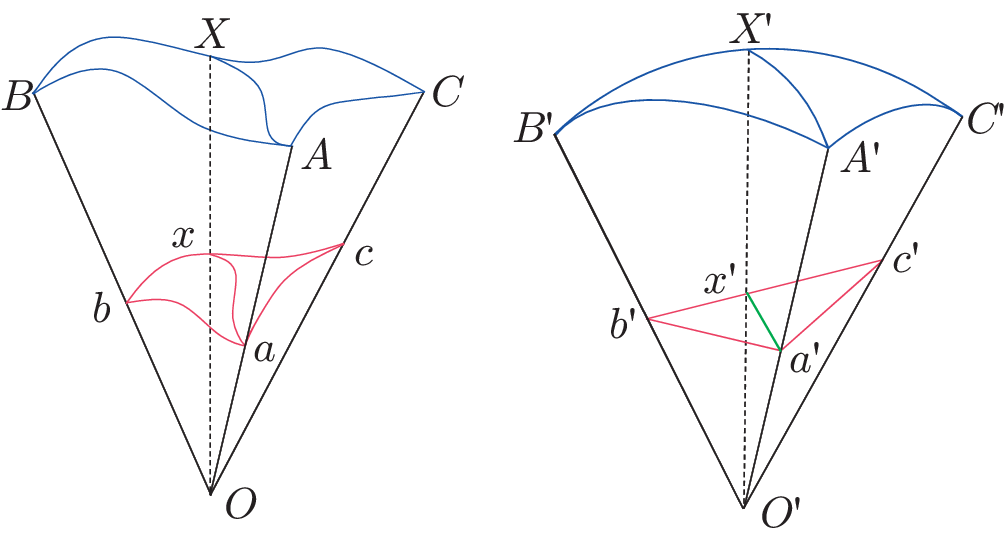}
\caption{The cone spanned by geodesics $\overline{ab}$, $\overline{ac}$, $\overline{bc}$ and
$\overline{ax}$ (left) and the cone spanned by a comparison triangle $\Delta a'b'c'$ (right).}
\label{fig:metric_cones}
\end{center}
\end{figure}

(1) Denote the mapped points of $a,b,c$ and $x$
by the projection $\tilde{\mathcal{X}}_\beta \rightarrow \mathcal{X}_\beta$ as $A,B,C$ and $X$, respectively,
as shown in Figure \ref{fig:metric_cones} (left).
Denote the origin of the metric cone as $O$.
If the sum of the lengths of the geodesics $\widetilde{AB}$, $\widetilde{AC}$ and $\widetilde{BC}$
in $\mathcal{X}_\beta$
exceeds $2\pi\beta$, it is easy to see that the cone spanned by $\overline{ab} \cup \overline{ac} \cup \overline{bc}$ becomes
CAT(0) and $\Delta abc$ satisfies the CAT(0) property.
Therefore, assume that $|\widetilde{AB}|+|\widetilde{AC}|+|\widetilde{BC}|\leq 2 \pi\beta$.

Next, let $\Delta a'b'c'$ be a comparison triangle of $\Delta abc$ and let $x'$ be a point on a geodesic $\overline{b'c'}$
such that $|\overline{bx}|=|\overline{b'x'}|$.
Thus, $|\overline{a'x'}|< |\overline{ax}|$.
Arrange the points $a'$, $b'$ and $c'$ in a three-dimensional Euclidean space with origin $O'$
such that the lengths of $\overline{O'a'}$, $\overline{O'b'}$ and $\overline{O'c'}$
are equal to the lengths of $\overline{Oa}$, $\overline{Ob}$ and $\overline{Oc}$, respectively.
Denote the radial projection of $a'$, $b'$, $c'$ and $x'$ to a unit sphere as $A'$, $B'$, $C'$ and $X'$, respectively,
as shown in Figure \ref{fig:metric_cones} (right).
By the definition of a metric cone, $|\overline{Ox}|=|\overline{O'x'}|$ and
the geodesics $\widetilde{A'B'}$, $\widetilde{A'C'}$, $\widetilde{B'C'}$ and $\widetilde{A'X'}$ in the unit sphere
are arcs satisfying
$|\widetilde{A'B'}|=|\widetilde{AB}|$, $|\widetilde{A'C'}|=|\widetilde{AC}|$, $|\widetilde{B'C'}|=|\widetilde{BC}|$
and $|\widetilde{B'X'}|=|\widetilde{BX}|$.

From the argument above, $|\widetilde{A'B'}|+|\widetilde{A'C'}|+|\widetilde{B'C'}|=
|\widetilde{AB}|+|\widetilde{AC}|+|\widetilde{BC}|\leq 2 \pi$.
Since the unit sphere has a positive constant curvature and $\mathcal{X}_\beta$ is CAT(0),
$|\widetilde{A'X'}|> |\widetilde{AX}|$.
However, since $|\overline{Oa}|=|\overline{O'a'}|$ and $|\overline{Ox}|=|\overline{O'x'}|$,
$|\widetilde{A'X'}|> |\widetilde{AX}|$ implies that $|\widetilde{a'x'}|> |\widetilde{ax}|$
by the property of a metric cone.
Thus, $\Delta abc$ has CAT(0) property and (1) of the theorem is proved.

(2) Assume that $0<\beta_1<\beta_2<\infty$ and a metric cone $\tilde{\mathcal{X}}_{\beta_1}$ is not CAT(0) for
proving the latter half of the theorem by contradiction.
Then, there is a geodesic triangle $\Delta a_1 b_1 c_1$ in $\tilde{\mathcal{X}}_{\beta_1}$ and
a point $x_1$ on the geodesic $\overline{b_1 c_1}$ such that the geodesic $\overline{a_1 x_1}$
is longer than the corresponding geodesic of a comparison triangle.
By defining $A_1, B_1, C_1, X_1,a'_1,b'_1,c'_1,x'_1,A'_1, B'_1, C'_1$ and $X'_1$ as above,
we can say that $|\widetilde{A_1 X_1}|>|\widetilde{A'_1 X'_1}|$.

Next, each of $A_1, B_1, C_1$ and $X_1$ corresponds to a point in $\mathcal{X}_{\beta_1}$ and
we can consider the corresponding points $A_2, B_2, C_2$ and $X_2$ in the other metric cone
$\tilde{\mathcal{X}}_{\beta_2}$.
When restricted to $\mathcal{X}_{\beta_1}$,
a geodesic $\widetilde{A_1 X_1}$ is just a rescaling of $\widetilde{A_2 X_2}$
and $|\widetilde{A_1 X_1}|=\frac{\beta_1}{\beta_2}|\widetilde{A_2 X_2}|$.

Now, $\Delta A_2'B_2'C_2'$ is a geodesic triangle on the unit sphere,
but after rescaling by $\frac{\beta_2}{\beta_1}$, we can get a geodesic triangle
$\Delta A_2''B_2''C_2''$ on a sphere of radius $\frac{\beta_2}{\beta_1}$
whose edges have the same length as $\Delta A_1'B_1'C_1'$.
By a known result on spherical triangles with the same edge lengths on different spheres,
a larger radius implies a ``thinner'' triangle and $|\widetilde{A_1 X_1}|<|\widetilde{A_2'' X_2''}|$
where $X_2''$ is a point on the geodesic $\widetilde{B_2'' C_2''}$ such that $|\widetilde{B_1 X_1}|=|\widetilde{B_2'' X_2''}|$.

Combining all the arguments gives
\begin{equation*}\textstyle |\widetilde{A_2 X_2}|=\frac{\beta_2}{\beta_1}|\widetilde{A_1 X_1}|>\frac{\beta_2}{\beta_1}|\widetilde{A'_1 X'_1}|
>\frac{\beta_2}{\beta_1}|\widetilde{A_2'' X_2''}|=|\widetilde{A'_2 X'_2}|.\end{equation*}
Select a non-degenerate geodesic triangle in $\tilde{\mathcal{X}}_{\beta_2}$ by selecting
arbitrary points $a_2,b_2$ and $c_2$ on the geodesics $\overline{OA_2}$, $\overline{OB_2}$ and $\overline{OC_2}$
in $\tilde{\mathcal{X}}_{\beta_2}$, respectively, and let $x_2$ be
the intersection point of $\overline{OX_2}$ and $\overline{b_2 c_2}$.
Then, by $|\widetilde{A_2 X_2}|>|\widetilde{A'_2 X'_2}|$, we can say that $|\overline{a_2 x_2}|>|\overline{a'_2 x'_2}|$.
This implies that $\tilde{\mathcal{X}}_{\beta_2}$ is not CAT(0) and (2) of the theorem is proved.

(3) For $k=0$, the statement holds by (1).
For $k>0$ and $\beta\leq\pi$, it is sufficient to prove for $\beta=\pi/\sqrt{k}$ by (2).
Let $\Delta abc$ be a geodesic triangle in $\tilde{\mathcal{X}}_\beta$ and
let $\Delta ABC$ be a geodesic triangle in $\mathcal{X}_\beta$.
Let $A,B,C$ be the projection of $a,b,c$, respectively.
If the perimeter of $\Delta ABC$ is longer than or equal to $2\pi$,
the cone spanned by the perimeter becomes CAT(0) by the same argument as that for (1).
Therefore, $\Delta abc$ is CAT(0) and satisfies the CAT(0) property.

If the perimeter of $\Delta ABC$ is smaller than $2\pi$, since $\mathcal{X}$ is CAT($k$) and
$\mathcal{X}_\beta$ is CAT(1),
for any $X\in \widetilde{BC}$, $\widetilde{BX}$ is shorter than the corresponding
great arc $\widetilde{B'X'}$ of a comparison triangle $\Delta A'B'C'$,
which is a spherical triangle on the unit sphere.
Since a comparison triangle $\Delta a'b'c'$ of $\Delta abc$ can be embedded
on the cone spanned by $\Delta A'B'C'$,
$\widetilde{bx}$ is shorter than the corresponding line segment $\widetilde{b'x'}$.
This means the $\Delta abc$ satisfies the CAT(0) property.
\qed

 \section{Proof of Theorem \ref{diam1}} \label{proof-diameter}

 (1) Although this is a known result, for example \cite{kendall-1990} \cite{espnola-fernandez-2009}, we show a short proof.
 Since a comparison triangle for the CAT($k$) property is on a sphere of radius $1/\sqrt{k}$,
 first consider the unit sphere $S^2$ and the geodesic distance $d$ on it.
 Take three points $a,b,c \in S^2$ and think of the convexity of $d(x,a)$ for $x\in \widetilde{bc}$.
 Without losing generality, assume that $a$ is on the plane $y=0$ and $\widetilde{bc}$ is on the plane
 $z=0$ and let $b=(\cos \theta_0, \sin \theta_0,0)$, $x=(\cos \theta, \sin \theta, 0)$ and
 $a=(\cos \psi, 0 ,\sin \psi)$ for $\theta_0, \theta\in (-\pi,\pi]$,$\psi\in[-\pi/2,\pi/2]$.

 Thus, $d(a,x)=\arccos(a^\top x)=\arccos(\cos \theta \cos\psi) $ and $d(b,x)=|\theta-\theta_0|$.
 Note that for $\psi=0$, $d(a,x)=|\theta|$ for $\psi=0$ is a convex of $d(b,x)$.
 For a $\psi \neq 0$, $d(a,x)\leq \pi/2$ for $x\in \widetilde{bc}$, $d(a,x)$ is a convex of $d(b,x)$
 iff $a^\top x=\cos \theta \cos \psi\geq 0$ since
 $\frac{\partial^2}{\partial \theta^2} d(a,x)= \cos(\theta) \cos(\psi) \sin^2(\psi) (1-\cos^2 \theta
 \cos^2\psi)^{-2/3}$.
 This means that if $d(a,x)\leq \pi/2$ for $x\in \widetilde{bc}$, $d(a,x)$ is a convex of $d(b,x)$.

 If $\mathcal{X}$ is CAT($k$) and has a diameter of at most $\pi/(2\sqrt{k})$,
 there is a comparison triangle $\Delta a'b'c'$ on a sphere of radius $1/\sqrt{k}$
 such that its perimeter is at most $3\pi/(2\sqrt{k})$ and $d(a',x')$ is a convex of $d(b',x')$
 for each $x'\in \widetilde{b'c'}$ because of the argument above after scaling by $1/\sqrt{k}$.

 (2) is well known. See \cite{espnola-fernandez-2009}.

 (3) We show an example of the probability measure with a three-point support on $S^2$
 such that the diameter is larger than $\pi/2$ but can be arbitrarily close to $\pi/2$
 and the uniqueness of the intrinsic $L_1$-mean fails.

 \begin{figure}[tb]
 \begin{center}
 \includegraphics[height=4.5cm]{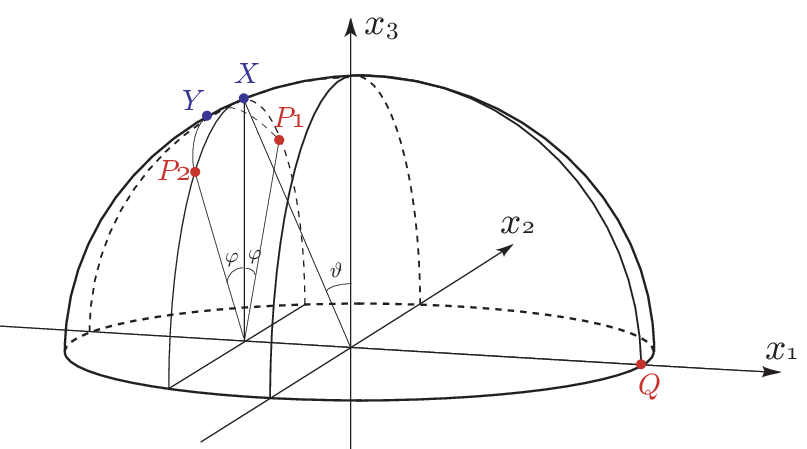}
 \caption{$P_1$,$P_2$,$Q$, $X$ and $Y$ on the unit hemisphere.}
 \label{fig:cat_hemisphere1}
 \end{center}
 \end{figure}

 Take $Q=(1,0,0)$, $P_1=(-\sin\theta,\cos\theta\sin\psi,\cos\theta \cos\psi)$,
 $P_2=(-\sin\theta,\cos\theta\sin\psi,\cos\theta \cos\psi)$ and $X=(-\sin\theta,0,\cos\theta)$
 with $\theta,\psi\in (0,\pi/2)$, as in Figure \ref{fig:cat_hemisphere1}.
 Let $Y=(-\sin\theta',0,\cos\theta')$ be the mid point of $\widetilde{P_1P_2}$
 for $\theta'>\theta$.
 Put the point masses $m_1$ at $Q$ and $M$ at $P_1$ and $P_2$, and
 assume that there is a unique intrinsic median $\mu$.

 By the symmetry, $\mu$ must be on the arc $\widetilde{QY}$, and if we change the ratio $M/m$,
 $\mu$ moves continuously on $\widetilde{QY}$.
 Thus, we can set $\mu=X$ by tuning $M/m$ adequately.
 However, $L_1$-dispersion from $X$ becomes
 $S_X=m d(\theta,X) +M d(P_1,X) + M d(P_2,X) = m(\pi/2+\theta)+2M\psi$, and
 $L_1$-dispersion from $P_1$ becomes
 $S_{P_1}=m d(\theta,P_1) +M d(P_2,P_1) = m(\pi/2+\theta)+2M\psi$.
 This contradicts the assumption of $X$ being the unique $L_1$-intrinsic mean.
 Since we can set $\theta$ and $\psi$ as arbitrarily small positive numbers, $D_{L_1}\leq  \pi/(2\sqrt{k})$. 
 However, by (1), $D_{L_1}\geq D_{\rm convex}\geq \pi/(2\sqrt{k})$;
 thus, $D_{L_1}=\pi/(2\sqrt{k})$.\qed


 \section{An upper bound of $D_{L_\gamma}$} \label{proof-gamma}

 \begin{thm}
 If $\mathcal{X}$ is a surface with a constant curvature $k>0$,
 \begin{equation*}D_{L_\gamma}\leq \frac{1}{\sqrt{k}}\left(\theta_0(\gamma)+\frac{\pi}{2}\right)\end{equation*}
 where $\theta_0(\gamma)$ is the inverse function of
 \begin{equation*}\gamma_0(\theta)=
 \left\{ \begin{array}{ll}
     \left[ \log_2 \frac{(\pi+2\theta) }{2\arccos\left\{\sin^2\theta + \cos^2\theta
     \left(1-1/2(1-\sin\theta)\right)^{1/2}\right\} }\right]^{-1}
  & \mbox{~for~} 0 \leq \theta \leq \pi/6,\\
     \left[\log_2\frac{\pi-2\theta}{\arccos(\sin^2\theta)}\right]^{-1}
  & \mbox{~for~} \pi/6< \theta \leq \pi/2
   \end{array} \right.
 \end{equation*}
 for $1\leq \gamma <2$ and $\theta_0(\gamma)=\pi/2$ for $\gamma\geq 2$.
 \end{thm}
 The graph of $\theta_0(\gamma)$ is shown in Figure \ref{fig:theta_0_graph}.

\vspace{0.5cm}
 \proof
 The proof is similar to that of Theorem \ref{thm:diameter}(3).
 We consider two cases of arrangement of three points $Q,P_1$ and $P_2$.

 C1:~$Q=(1,0,0)$, $P_1=(-\sin\theta,\cos\theta\sin\psi,\cos\theta\cos\psi)$,\\
 ~~~~~~~~$P_2=(-\sin\theta,-\cos\theta\sin\psi,\cos\theta\cos\psi)$ where
 $\psi=\arccos \left\{\left(1-\frac{1}{2(1-\sin\theta)}\right)^{1/2}\right\}$, as
 shown in Figure \ref{fig:cat_hemisphere1}.
 This satisfies $|\widetilde{P_1Q}|=|\widetilde{P_1P_2}|$.
 \\

\begin{figure}[tb]
 \begin{center}
 \includegraphics[height=4cm]{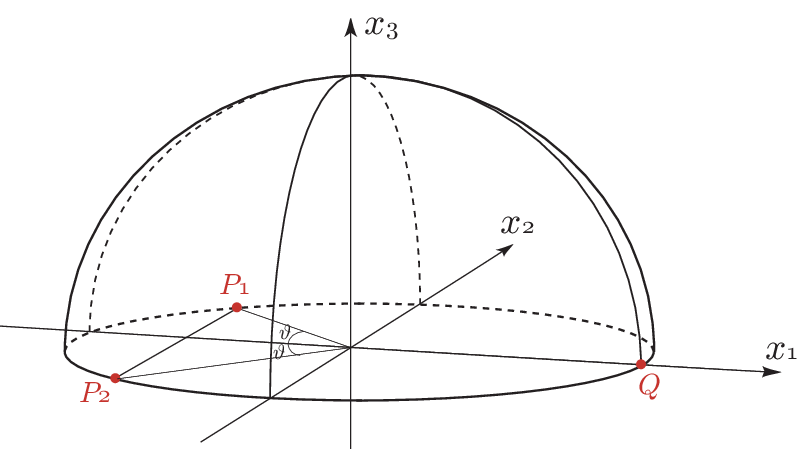}
 \caption{$P_1$,$P_2$ and $Q$ on the unit hemisphere.}
 \label{fig:cat_hemisphere2}
 \end{center}
 \end{figure}

 C2: $Q=(1,0,0)$, $P_1=(-\sin\theta,\cos\theta,0)$, and $P_2=(-\sin\theta,-\cos\theta,0)$,
 as shown in Figure \ref{fig:cat_hemisphere2}.
 \\

 We put point masses $m$ at $Q$ and $M$ at $P_1$ and $P_2$.

 For $1\leq \gamma \leq \left\{\log_2\left(\frac{2\pi}{3\arccos{(1/4)}}\right)\right\}^{-1}$, we consider
 C1.
 As in the proof of Theorem \ref{thm:diameter}(3), we can set $\mu=X=(-\sin\theta,0,\cos\theta)$.
 Let $S_X$ and $S_{P_1}$ denote $L_\gamma$-dispersion from $X$ and $P_1$, respectively.
 Then,
 \begin{eqnarray*}
 S_{X}&= m(\pi/2+\theta) ^\gamma + 2M\{\arccos(1-\cos^2\theta(1-\cos\psi))\},\\
 S_{P_1}&=  m(\pi/2+\theta) ^\gamma + M\{\arccos(1-2\cos^2\theta \sin^2\psi)\}.
 \end{eqnarray*}
 Therefore, $S_X< S_{P_1}$ is equivalent to
 \begin{equation*}\arccos(1-2\cos^2\theta\sin^2\psi)< 2^{1/\gamma} \arccos\{1-\cos^2\theta(1-\cos\psi)\}.\end{equation*}

 By setting $\psi=\arccos\{(1-1/2(1-\sin\theta))^{1/2}\}$, this is equivalent to
 $\gamma< \gamma_0(\theta)$ and also $\theta < \theta_0(\gamma)$.
 Thus, if we set $\theta\geq \theta_0(\gamma)$, C1 becomes an example of a non-unique intrinsic
 $L_1$-mean of diameter $\theta+\pi/2$.

 For C2, $S_X< S_{P_1}$ is equivalent to
 $\pi-2\theta < 2^{1/\gamma} \arccos(\sin^2\theta)$, and we can
 prove that it becomes a similar example.
 After scaling by $1/\sqrt{k}$, these examples give the upper bound on $D_{L_\gamma}$. \qed

\section*{Acknowledgements}
Funding was provided by JST, PRESTO
(JPMJPR14E3), JSPS, KAKENHI (26280009,16K02843) and RIKEN, AIP Japan.
The first author would like to thank Masayuki Sakai, Takaaki Koike and
Tatsuhiro Aoshima for their excellent computation and visualization of
the results. He also appreciates Reiko Miyaoka and Hiroshi Kokubu for
their helpful and encouraging advice.

%
%


%

\bibliographystyle{spbasic}      
\bibliography{Kobayashi_Wynn_SC}   

\end{document}